\ifx\shlhetal\undefinedcontrolsequence\let\shlhetal\relax\fi

\input amstex
 \loadmsbm
\NoBlackBoxes
\documentstyle {amsppt}
\topmatter
\title {{\it On Monk's Questions}\\
Sh479} \endtitle
\author {Saharon Shelah \thanks {\null\newline
Partially supported by the Deutsche
Forschungsgemeinschaft, Grant Ko 490/7-1. \null\newline
I would like to thank Alice Leonhardt for the beautiful typing. \null\newline
Publication 479} \endthanks } \endauthor
\affil {Institute of Mathematics\\
The Hebrew University\\
Jerusalem, Israel
\medskip
Rutgers University\\
Department of Mathematics\\
New Brunswick, NJ USA} \endaffil
\endtopmatter
\document
%
%
\def\renewcommand{\newcommand}	       
\edef\cite{\the\catcode`@}%
\catcode`@ = 11
\let\@oldatcatcode = \cite
\chardef\@letter = 11
\chardef\@other = 12
%
%
%
%
\def\@innerdef#1#2{\edef#1{\expandafter\noexpand\csname #2\endcsname}}%
%
%
\@innerdef\@innernewcount{newcount}%
\@innerdef\@innernewdimen{newdimen}%
\@innerdef\@innernewif{newif}%
\@innerdef\@innernewwrite{newwrite}%
%
%
%
\def\@gobble#1{}%
%
%
%
\ifx\inputlineno\@undefined
   \let\@linenumber = \empty 
\else
   \def\@linenumber{\the\inputlineno:\space}%
\fi
%
%
%
\def\@futurenonspacelet#1{\def\cs{#1}%
   \afterassignment\@stepone\let\@nexttoken=
}%
\begingroup 
\def\\{\global\let\@stoken= }%
\\ 
\endgroup
\def\@stepone{\expandafter\futurelet\cs\@steptwo}%
\def\@steptwo{\expandafter\ifx\cs\@stoken\let\@@next=\@stepthree
   \else\let\@@next=\@nexttoken\fi \@@next}%
\def\@stepthree{\afterassignment\@stepone\let\@@next= }%
%
%
%
\def\@getoptionalarg#1{%
   \let\@optionaltemp = #1%
   \let\@optionalnext = \relax
   \@futurenonspacelet\@optionalnext\@bracketcheck
}%
%
%
\def\@bracketcheck{%
   \ifx [\@optionalnext
      \expandafter\@@getoptionalarg
   \else
      \let\@optionalarg = \empty
      \expandafter\@optionaltemp
   \fi
}%
\def\@@getoptionalarg[#1]{%
   \def\@optionalarg{#1}%
   \@optionaltemp
}%
%
%
%
\def\@nnil{\@nil}%
\def\@fornoop#1\@@#2#3{}%
\def\@for#1:=#2\do#3{%
   \edef\@fortmp{#2}%
   \ifx\@fortmp\empty \else
      \expandafter\@forloop#2,\@nil,\@nil\@@#1{#3}%
   \fi
}%
\def\@forloop#1,#2,#3\@@#4#5{\def#4{#1}\ifx #4\@nnil \else
       #5\def#4{#2}\ifx #4\@nnil \else#5\@iforloop #3\@@#4{#5}\fi\fi
}%
\def\@iforloop#1,#2\@@#3#4{\def#3{#1}\ifx #3\@nnil
       \let\@nextwhile=\@fornoop \else
      #4\relax\let\@nextwhile=\@iforloop\fi\@nextwhile#2\@@#3{#4}%
}%
%
%
%
\@innernewif\if@fileexists
\def\@testfileexistence{\@getoptionalarg\@finishtestfileexistence}%
\def\@finishtestfileexistence#1{%
   \begingroup
      \def\extension{#1}%
      \immediate\openin0 =
         \ifx\@optionalarg\empty\jobname\else\@optionalarg\fi
         \ifx\extension\empty \else .#1\fi
         \space
      \ifeof 0
         \global\@fileexistsfalse
      \else
         \global\@fileexiststrue
      \fi
      \immediate\closein0
   \endgroup
}%
%
%
%
%
\def\bibliographystyle#1{%
   \@readauxfile
   \@writeaux{\string\bibstyle{#1}}%
}%
\let\bibstyle = \@gobble
%
%
\let\bblfilebasename = \jobname
\def\bibliography#1{%
   \@readauxfile
   \@writeaux{\string\bibdata{#1}}%
   \@testfileexistence[\bblfilebasename]{bbl}%
   \if@fileexists
      \nobreak
      \@readbblfile
   \fi
}%
\let\bibdata = \@gobble
%
%
\def\nocite#1{%
   \@readauxfile
   \@writeaux{\string\citation{#1}}%
}%
\@innernewif\if@notfirstcitation
%
%
\def\cite{\@getoptionalarg\@cite}%
%
%
\def\@cite#1{%
   \let\@citenotetext = \@optionalarg
   \printcitestart
   \nocite{#1}%
   \@notfirstcitationfalse
   \@for \@citation :=#1\do
   {%
      \expandafter\@onecitation\@citation\@@
   }%
   \ifx\empty\@citenotetext\else
      \printcitenote{\@citenotetext}%
   \fi
   \printcitefinish
}%
\def\@onecitation#1\@@{%
   \if@notfirstcitation
      \printbetweencitations
   \fi
   \expandafter \ifx \csname\@citelabel{#1}\endcsname \relax
      \if@citewarning
         \message{\@linenumber Undefined citation `#1'.}%
      \fi
      \expandafter\gdef\csname\@citelabel{#1}\endcsname{%
         {\tt
            \escapechar = -1
            \nobreak\hskip0pt
            \expandafter\string\csname#1\endcsname
            \nobreak\hskip0pt
         }%
      }%
   \fi
   \csname\@citelabel{#1}\endcsname
   \@notfirstcitationtrue
}%
%
%
\def\@citelabel#1{b@#1}%
%
%
\def\@citedef#1#2{\expandafter\gdef\csname\@citelabel{#1}\endcsname{#2}}%
%
%
%
\def\@readbblfile{%
   \ifx\@itemnum\@undefined
      \@innernewcount\@itemnum
   \fi
   \begingroup
      \def\begin##1##2{%
         \setbox0 = \hbox{\biblabelcontents{##2}}%
         \biblabelwidth = \wd0
      }%
      \def\end##1{}
      %
      %
      \@itemnum = 0
      \def\bibitem{\@getoptionalarg\@bibitem}%
      \def\@bibitem{%
         \ifx\@optionalarg\empty
            \expandafter\@numberedbibitem
         \else
            \expandafter\@alphabibitem
         \fi
      }%
      \def\@alphabibitem##1{%
         \expandafter \xdef\csname\@citelabel{##1}\endcsname {\@optionalarg}%
         \ifx\biblabelprecontents\@undefined
            \let\biblabelprecontents = \relax
         \fi
         \ifx\biblabelpostcontents\@undefined
            \let\biblabelpostcontents = \hss
         \fi
         \@finishbibitem{##1}%
      }%
      \def\@numberedbibitem##1{%
         \advance\@itemnum by 1
         \expandafter \xdef\csname\@citelabel{##1}\endcsname{\number\@itemnum}%
         \ifx\biblabelprecontents\@undefined
            \let\biblabelprecontents = \hss
         \fi
         \ifx\biblabelpostcontents\@undefined
            \let\biblabelpostcontents = \relax
         \fi
         \@finishbibitem{##1}%
      }%
      \def\@finishbibitem##1{%
         \biblabelprint{\csname\@citelabel{##1}\endcsname}%
         \@writeaux{\string\@citedef{##1}{\csname\@citelabel{##1}\endcsname}}%
         \ignorespaces
      }%
      %
      %
      \let\em = \bblem
      \let\newblock = \bblnewblock
      \let\sc = \bblsc
      \frenchspacing
      \clubpenalty = 4000 \widowpenalty = 4000
      \tolerance = 10000 \hfuzz = .5pt
      \everypar = {\hangindent = \biblabelwidth
                      \advance\hangindent by \biblabelextraspace}%
      \bblrm
      \parskip = 1.5ex plus .5ex minus .5ex
      \biblabelextraspace = .5em
      \bblhook
      \input \bblfilebasename.bbl
    \endgroup
}%
%
%
\@innernewdimen\biblabelwidth
\@innernewdimen\biblabelextraspace
%
%
%
\def\biblabelprint#1{%
   \noindent
   \hbox to \biblabelwidth{%
      \biblabelprecontents
      \biblabelcontents{#1}%
      \biblabelpostcontents
   }%
   \kern\biblabelextraspace
}%
%
%
%
\def\biblabelcontents#1{{\bblrm [#1]}}%
%
%
\def\bblrm{\rm}%
%
%
\def\bblem{\it}%
%
%
\def\bblsc{\ifx\@scfont\@undefined
              \font\@scfont = cmcsc10
           \fi
           \@scfont
}%
%
%
\def\bblnewblock{\hskip .11em plus .33em minus .07em }%
%
%
\let\bblhook = \empty
%
%
%
\def\printcitestart{[}
\def\printcitefinish{]}
\def\printbetweencitations{, }
\def\printcitenote#1{, #1}
%
%
%
\let\citation = \@gobble
%
%
%
\@innernewcount\@numparams
%
%
\def\newcommand#1{%
   \def\@commandname{#1}%
   \@getoptionalarg\@continuenewcommand
}%
%
%
\def\@continuenewcommand{%
   \@numparams = \ifx\@optionalarg\empty 0\else\@optionalarg \fi \relax
   \@newcommand
}%
%
%
\def\@newcommand#1{%
   \def\@startdef{\expandafter\edef\@commandname}%
   \ifnum\@numparams=0
      \let\@paramdef = \empty
   \else
      \ifnum\@numparams>9
         \errmessage{\the\@numparams\space is too many parameters}%
      \else
         \ifnum\@numparams<0
            \errmessage{\the\@numparams\space is too few parameters}%
         \else
            \edef\@paramdef{%
               \ifcase\@numparams
                  \empty  No arguments.
               \or ####1%
               \or ####1####2%
               \or ####1####2####3%
               \or ####1####2####3####4%
               \or ####1####2####3####4####5%
               \or ####1####2####3####4####5####6%
               \or ####1####2####3####4####5####6####7%
               \or ####1####2####3####4####5####6####7####8%
               \or ####1####2####3####4####5####6####7####8####9%
               \fi
            }%
         \fi
      \fi
   \fi
   \expandafter\@startdef\@paramdef{#1}%
}%
%
%
%
%
\def\@readauxfile{%
   \if@auxfiledone \else 
      \global\@auxfiledonetrue
      \@testfileexistence{aux}%
      \if@fileexists
         \begingroup
            \endlinechar = -1
            \catcode`@ = 11
            \input \jobname.aux
         \endgroup
      \else
         \message{\@undefinedmessage}%
         \global\@citewarningfalse
      \fi
      \immediate\openout\@auxfile = \jobname.aux
   \fi
}%
%
%
\newif\if@auxfiledone
\ifx\noauxfile\@undefined \else \@auxfiledonetrue\fi
%
%
%
%
\@innernewwrite\@auxfile
\def\@writeaux#1{\ifx\noauxfile\@undefined \write\@auxfile{#1}\fi}%
%
%
%
\ifx\@undefinedmessage\@undefined
   \def\@undefinedmessage{No .aux file; I won't give you warnings about
                          undefined citations.}%
\fi
%
%
\@innernewif\if@citewarning
\ifx\noauxfile\@undefined \@citewarningtrue\fi
%
%
%
\catcode`@ = \@oldatcatcode

%
%

\newpage
\head{Anotated Content} \endhead
\bigskip

\noindent
\S1  Introduction
\bigskip

\noindent
\S2  Existence of subalgebras with a preassigned algebraic density
\smallskip

[We first note (in 2.1) that if $\pi(B) \ge \theta = \text{cf}(\theta)$
then for some $B' \subseteq B$ we have $\pi(B') = \theta$.  Call this 
statement $(*)$.  Then we give a
criterion for \newline
$\pi(B) = \mu > \text{ cf}(\mu)$ (in 2.2) and conclude for
singular $\mu$ that for a club of $\theta < \mu$ the $(*)$ above 
holds (2.2A), and
investigate the criterion (in 2.3).  Our main aim is, starting with
$\mu = \mu^{< \mu},\text{cf}(\lambda) < \lambda$, to force the existence
of a Boolean algebra $B$ with $\pi(B) > \theta$ but for no $B' \subseteq B$
do we have $\pi(B') = \lambda$ (in fact $(\exists B' \subseteq B)
[\pi(B') = \theta \Leftrightarrow \theta = \text{ cf}(\theta) \vee
\text{ cf}(\theta) \le \mu]$ for every $\theta \le |B|$.  Toward this, we
define the forcing (Definition 2.5: how $\langle x_\alpha:\alpha \in W^P 
\rangle$ generate a Boolean algebra, $BA[p],W^p \in [\lambda]^{< \mu}$ with
$x_\alpha > \theta$ has no non-zero member of $\langle x_\beta:\beta \in
W^p \cap \alpha \rangle_{BA[p]}$ below it).  We prove the expected properties
of the generic (2.6), also the forcing has the expected properties 
($\mu$-complete, $\mu^+$-c.c) (in 2.7).  The main theorem (2.9) stated, 
the main point brings that
if $\mu < \text{ cf}(\theta) < \theta$ for $B \subseteq BA[G],
\pi(B) \ne \theta$; we use the criterion from above, a lemma related to
$\triangle$-systems (see \cite{Sh:430},6.6D, \cite{Sh:513},6.1) quoted in
2.4; to reduce the problem to some special amalgamation of finitely many
copies (the exact number is in relation to the arity of the term defining 
the relevant
elements from the $x_\alpha$'s).  The existence of such amalgamation was
done separately earlier (2.8).  

Lastly in 2.10 we show that the cf$(\theta) > \mu$ above was necessary by
proving the existence of a subalgebra with prescribed singular algebraic 
density
$\lambda:\pi(B) > \lambda$ and $(\forall \mu < \lambda)[\mu^{< \text{ cf}
(\lambda)} < \lambda]$].
\bigskip

\noindent
\S3 On $\pi$ and $\pi \chi$ of Products of Boolean Algebras
\smallskip

[If e.g. $\aleph_0 < \kappa = \text{ cf}(\chi) < \chi < \lambda =
\text{ cf}(\lambda) < \chi^\kappa,(\forall \theta < \chi)(\theta^\kappa
< \chi)$ we show that for some Boolean algebras $B_i$ (for $i < \kappa):
\chi - \dsize \sum_{i < \chi} \pi \chi(B_i) < \lambda$ but (for $D$ a regular
ultrafilter on $\kappa$) $\lambda = \pi \chi (\dsize \prod_{i < \kappa}
B_i/D)$ but $\dsize \prod_{i < \kappa}(\pi \chi(B_i))/D = \chi^\kappa$.  For
this we use interval Boolean algebras on order of the form $\lambda_i \times
\Bbb Q$.

We also prove for infinite Boolean algebras $B_i$ (for $i < \kappa$) and $D$
an ultrafilter on $\kappa$, if $n_i < \aleph_0,\mu = \dsize \prod_{i < \kappa}
n_i/D$ is regular (infinite) cardinal then $\pi \chi(\dsize \prod_{i < \kappa}
B_i/D) \ge \mu$.
\newpage

\head {\S1 Introduction}\endhead
\medskip
\noindent  Monk \cite{M} asks: (problems 13, 15 in his list; $\pi$ is the
algebraic density, see 1.1 below) \newline

For a (Boolean algebra) $B$, $\aleph_0 \le \theta \le \pi(B)$, does $B$ have
a subalgebra $B^\prime$ with $\pi(B^\prime) = \theta$? \newline
\medskip

If $\theta$ is regular the answer is easily positive (see 2.1), we show that
in general it may be negative (see 2.9(3)), but for quite many singular
cardinals - it is positive (2.10); the theorems are quite complementary.  
This is dealt with in $\S 2$.\newline

In $\S 3$ we mainly deal with $\pi \chi$ (see Definition 3.2) show that the
$\pi \chi$ of an ultraproduct of Boolean algebras is not necessarily the
ultraproduct of the $\pi \chi$'s.  Note: in Koppelberg Shelah \cite{KpSh:415}
,Theorem 1.1 we prove that if SCH holds,
$\pi(B_i) > 2^\kappa$ for $i < \kappa$ then 
$\pi \left( \dsize \prod_{i<\kappa}
B_i/D \right) = \dsize \prod_{i<\kappa}(\pi(B_i))/D$. \newline
We also prove that for infinite Boolean algebras $A_i(i < \kappa)$ and a
non-principal ultrafilter $D$ on $\kappa$: if $n_i < \aleph_0$ for 
$i < \kappa$ and $\mu =: \dsize \prod_{i < \kappa} n_i/D$ is regular, 
then $\pi\chi(A) \ge \mu$.  Here $A =: \dsize \prod_{i < \kappa} A_i/D$.
By a theorem of Peterson the regularity of $\mu$ is needed.
\bigskip

\noindent 1.1 \underbar{Notation:}  \, Boolean algebras are denoted by $B$ and
sometimes $A$.\newline
\noindent  For a Boolean algebra $B$
$$
B^+ =: \{ x \in B:x \ne 0 \}
$$
$$
\pi(B) =: \text{ Min } \{ |X|:X \subseteq B^+ \text{ is such that } \forall y
\in B^+ \exists x \in X[x \le y] \}.
$$
$X$ like that is called dense in $B$.  More generally if $X,Y \subseteq B$
we say $X$ is dense in $Y$ if $y \in Y \and y \ne 0 \Rightarrow (\exists
x \in X)(0 < x \le y)$.  For a $Y \subseteq B$, $\langle Y \rangle_B$
is the subalgebra of $B$ which
$Y$ generates.\newline
\noindent  $0_A$ is the constant function with domain $A$ and value
zero.\newline
\noindent  $1_A$ is defined similarly.
\newpage

\head {\S2  Existence of subalgebras with a preassigned algebraic
density} \endhead
\medskip

\noindent Note
\demo{2.1 Observation}  If $\pi(B) > \theta = \text{ cf}(\theta) \ge |Y| 
+ \aleph_0$ and $Y \subseteq B$ \underbar{then} for some
subalgebra $A$ of $B$, $Y \subseteq A$ and $\pi(A)
= \theta$.
\enddemo
\bigskip

\demo{Proof}  Without loss of generality $|Y| = \theta$.  Let
$Y = \{ y_\alpha:\alpha < \theta \}$.  Choose by induction
on $\alpha \le \theta$ subalgebras $A_\alpha$ of $B$, increasing continuous in
$\alpha$, $|A_\alpha| < \theta$, $y_\alpha \in A_{\alpha +1}$ such that: for
each $\alpha < \theta$, some $x_\alpha \in A^+_{\alpha +1}$ is not above
any $y \in A^+_\alpha$.  This is possible because for no $\alpha < \theta$ 
can $A^+_\alpha$ be dense in $B$.  \newline
Now $A = A_\theta$ is as required.\hfill$\square_{2.1}$
\enddemo
\bigskip
  
\proclaim{2.2 Claim}  Assume $B$ is a Boolean algebra, $\pi(B)=\mu >
\text{ cf}(\mu) \ge \aleph_0$ (see Definition~1.1).  \underbar{Then} 
for arbitrarily large regular $\theta<\mu$
\endproclaim

$$
\text{ for some set } Y \text{ we have:} \tag"($\ast$)$^B_\theta$"
$$

$$
Y\subseteq B^+, |Y| = \theta, \text{ and there is
\underbar{no} } Z\subseteq B^+ \tag"($\ast$)$^B_\theta [Y]$"
$$

$$
\text{ of cardinality } <\theta, \text{ dense in } Y
\left( \text{ i.e. } \bigwedge _{y\in Y} \bigvee
_{z\in Z}z\le y \right).
$$
\bigskip

\demo{2.2A Conclusion}  If $B$ is a Boolean algebra, $\pi(B) = \mu >
\text{ cf}(\mu) > \aleph_0$ and \newline
$\langle \mu_\zeta:\zeta < \text{ cf}(\mu)
\rangle$ is increasing continuously with limit $\mu$ (so $\mu_\zeta < \mu$)
\underbar{then} for some club $C$ of cf$(\mu)$ for every $\zeta \in C$ for
some $B' \subseteq B$ we have $\pi(B') = \mu_\zeta$.
\enddemo
\bigskip

\demo{Proof of 2.2} Let $Z^\ast \subseteq B^+$ be dense, $|Z^\ast|=\mu$.

If the conclusion fails, then for some $\theta^\ast < \mu$, for no regular
$\theta \in (\theta^\ast,\mu)$ does $(*)^B_\theta$ hold.  We now assume we
chose such $\theta^*$, and show by induction on $\lambda \le \mu$ that:

$$
\text{ if } Y\subseteq B^+,\, |Y|\le \lambda \text{ then for some } Z\subseteq
B^+, \, |Z|\le \theta^\ast \text{ and } Z \text{ is dense in } Y.
\tag"$\otimes_\lambda:$"
$$
\bigskip
\demo{First Case}  $\lambda \le \theta^\ast$. \enddemo

Let $Z = Y$.
\bigskip
\demo{Second Case}  $\theta^\ast < \lambda \le \mu$ and cf$(\lambda) < 
\lambda$.\newline

Let $Y = \bigcup\{Y_\zeta:\zeta < 
\text{ cf}(\lambda)\}$,\,$|Y_\zeta|<\lambda$. 
By the
induction hypothesis for each $\zeta < \text{ cf}(\lambda)$ there is
$Z_\zeta\subseteq B^+$ of cardinality $\le \theta^*$ which is
dense in $Y_\zeta$.

Now $Z' =: \dsize\bigcup_{\zeta < \text{ cf}(\lambda)}Z_\zeta$ has cardinality
$\le \theta^\ast + \text{ cf}(\lambda) < \lambda$, hence by the induction
hypothesis there is $Z \subseteq B^+$ dense in $Z'$ with $|Z| \le \theta^*$.
Easily $Z$ is dense in $Y$,\,$|Z| \le \theta^\ast$ and
$Z \subseteq B^+$ so we finish the case.
\enddemo  
\bigskip

\demo{Third Case}  $\theta^\ast < \lambda \le \mu$, \, $\lambda$ regular.
\newline

If for this $Y$, $(*)^B_\lambda[Y]$ holds,
we get the conclusion of the claim.  We are assuming not so; so there is 
$Z^\prime \subseteq B^+$, $|Z^\prime| < \lambda$, $Z'$ dense in $Y$.  Apply 
the induction hypothesis to $Z^\prime$ and get $Z$ as required. \newline
\noindent So we have proved $\otimes_\lambda$.
\medskip
We apply $\otimes_\lambda$ to $\lambda=\mu$, \, $Y=Z^\ast$ and 
get a contradiction. \hfill$\square_{2.2}$
\enddemo
\bigskip

\proclaim{2.3 Claim}  1) If $B$, $\mu$, $\theta$, $Y$ are as in 2.2
(so $(*)^B_\theta[Y]$ and $\theta$ is regular)  
\underbar{then} we can find $\bar y = \langle y_\alpha:\alpha < \theta
\rangle$ contained in $B^+$ and a proper $\theta$-complete filter 
$D$ on $\theta$ containing all cobounded subsets of $\theta$ such that:
$$
\text{ for every } z \in B^+ \text { , } \{ \alpha <
\theta:z \le y_\alpha \} = \emptyset \text{ mod } D.
\tag"$\otimes^B_{\bar y,D}$"
$$
2) If in addition $\theta$ is a successor cardinal \underbar{then} we can
demand that $D$ is normal.\endproclaim 
\bigskip

\remark{Remark}  Part (2) is for curiousity only.\endremark
\bigskip
\demo{Proof}  1) Let $Y = \{ y_\alpha:\alpha < \theta \}$.  Define $D$:
$$
\split
\text{ for } {\Cal U} \subseteq \theta \, ; \, {\Cal U} \in D
\text{ \underbar{iff} for some } 0 < \zeta < \theta \text{ and } z_\epsilon
\in B^+ \text{ for }\\
\qquad \epsilon < \zeta, \text{ we have } {\Cal U} \supseteq
\{ \alpha < \theta:\dsize\bigwedge_{\epsilon < \zeta}z_\epsilon \nleq
y_\alpha \}.
\endsplit
$$
Trivially $D$ is closed under supersets and intersections of $< \theta$
members and every cobounded subset of $\theta$ belongs to it.
Now $\emptyset \notin
D$ because $(*)^B_\theta[Y]$.\newline
\noindent 2)  Let $\theta = \sigma^+$.  Assume there are no such
$\bar y$, \, $D$.  We try to choose by induction on  $n < \omega$, \,
$Y^n_\alpha$ $(\alpha < \theta)$ and club $E_n$ of $\theta$ such that:
\roster
\item"{\bf(a)}" $Y^n_\alpha$ is a subset of $B^+$ of
cardinality $< \theta$, increasing continuous in $\alpha$
\item"{\bf(b)}" $Y^n_\alpha \subseteq Y^n_{\alpha +1}$
\item"{\bf(c)}" $Y^0_\alpha = \{ y_\beta:\beta < \alpha \}$ (taken from part
(1))
\item"{\bf(d)}" $E_n$ is a club of $\theta$, \, $E_{n+1} \subseteq E_n$, \,
$E_0 = \{ \delta < \theta:\delta \text{ divisible by } \sigma \}$
\item"{\bf(e)}" if $\delta \in E_{n+1}$ and $\delta \le \alpha < \text{ Min}
(E_n \backslash (\delta +1))$ \underbar{then} for every $y \in Y^n_\alpha$
there is $z \in Y^{n+1}_\delta$, \, $z \le y$.
\endroster

\noindent If we succeed, let $\beta^\ast = \dsize \bigcup_{n<\omega} \,
\text{ Min}(E_n)$ \, $(< \theta)$,
\,  and we shall prove that $\dsize \bigcup_{n<\omega} Y^n_{\beta^{\ast}}$
is dense in $Y$, getting a contradiction.  For every $y \in \dsize \bigcup
\Sb n< \omega \\ \alpha < \theta \endSb Y^n_\alpha$ \, let $\beta(y)$ be
the minimal $\beta < \theta$ such that $\left( \exists z \in \dsize \bigcup
\Sb n<\omega \endSb Y^n_\beta \right)(z \le y)$.  Now $\beta$ is well 
defined as $\langle \dsize \bigcup_{n<\omega}Y^n_\beta:
\beta < \theta \rangle$ is
increasing continuous and $y \in \dsize \bigcup_{\beta < \theta}
\dsize \bigcup_{n < \omega} Y^n_\beta$.
If $\beta(y) \le \beta^\ast$ for every \newline
$y \in Y \, \biggl( \subseteq
\dsize \bigcup_{\alpha < \theta} Y^0_\alpha \biggr)$ we are done, assume not,
so some $y^\ast \in Y = \dsize \bigcup_{\alpha < \theta} Y^0_\alpha$
exemplifies this.  Now let $\beta = \beta(y^\ast)$ and $z \in \dsize
\bigcup_{n < \omega} Y^n_\beta$ exemplifies this.  Clearly
$\langle \sup(\beta \cap E_n):n
< \omega \rangle$ is well defined; clearly it is a
non-increasing sequence of ordinals
hence eventually constant, say $n \ge n^\ast \Rightarrow \sup(\beta \cap
E_n) = \gamma$.  Now, without loss of generality $z \in Y^{n{^\ast}}_\beta$ 
(by clause (b)); note $\gamma \in E_n$ for $n \ge n^\ast$ (hence for every
$n$).  But by clause (e) there is $z^\prime \in Y^{n{^\ast}+1}_\gamma$,  $z^\prime \le z$, \,
contradicting the choice of $\beta$.

So we cannot carry the construction, so we are stuck at some $n$.
Fix such an $n$.  Let
$E_n \cup \{ 0 \} = \{ \delta_\epsilon:\epsilon < \theta \}$ (increasing with
$\epsilon$).  Let  $Y^n_{\delta _{\epsilon +1}}\backslash
Y^n_{\delta_\epsilon} \subseteq
\{ y^\epsilon _\zeta : \zeta < \sigma \}$.
For each $\zeta < \sigma$, let $D_\zeta$ be the normal filter generated by
the family of subsets of $\theta$ of the form $\{ \varepsilon < \theta:
z \nleq y^\varepsilon_\zeta \}$ for $z \in B^+$.  If for every
$\zeta < \sigma$, $\emptyset \in D_\zeta$ we can define $Y^{n+1}_\alpha$,
$E_{n+1}$, contradiction.  So for some $\zeta$, $\bar y^\zeta =: \langle
Y^\varepsilon_\zeta:\zeta < \sigma \rangle,D_\zeta$ are as required in
$\bigotimes^B_{\bar y^\zeta,D_\zeta}$.  \hfill$\square_{2.3}$
\enddemo
\bigskip

\proclaim{2.4 Claim}  Suppose $D$ is a filter on $\theta$,
$\sigma$-complete, $\theta = \text{ cf}(\theta) \ge \sigma > 2^\kappa$,
and for each $\alpha < \theta$,
$\bar \beta^\alpha = \langle \beta^\alpha_\epsilon: 
\epsilon < \kappa \rangle$ is a
sequence of ordinals. \underbar{Then} for every ${\Cal U} \subseteq
\theta$, ${\Cal U} \ne \emptyset$ mod $D$ there are $\langle \beta^\ast
_\epsilon:\epsilon < \kappa \rangle$ (a sequence of ordinals) and 
$w \subseteq \kappa$ such that:

\roster
\item "{\bf(a)}"  $\epsilon \in \kappa \backslash w \Rightarrow \text{ cf}
(\beta^\ast_\epsilon) \le \theta$
\item "{\bf(b)}"  if $\dsize \bigwedge_{\alpha < \sigma}|\alpha|^\kappa <
\sigma$ then: \newline
$\varepsilon \in \kappa \backslash w \Rightarrow \sigma \le \text{ cf}(\beta
^\ast_\varepsilon)$
\item "{\bf(c)}"  if $\beta^\prime_\epsilon \le \beta^\ast_\epsilon$ for
all $\varepsilon$ and
$[\varepsilon \in w \equiv \beta^\prime_{\epsilon} = \beta^\ast_{\epsilon}]$
then \newline
$\left\{ \alpha \in {\Cal U}:\beta^\prime_\epsilon \le \beta^\alpha_\epsilon
\le \beta^\ast_\epsilon \text{ for all } \varepsilon \text{ and } 
[\varepsilon \in w \equiv \beta^\alpha_\epsilon =
\beta^\ast_\epsilon] \right\} \ne \emptyset$ mod $D$.
\endroster
\endproclaim
\bigskip

\demo{Proof}  \cite{Sh:430},6.1D and better presented in \cite{Sh:513},6.1.
\enddemo
\bigskip

\definition{2.5 Definition}  1) If $F \subseteq {}^w2$ let \newline
$c \ell (F) = \{ g \in {}^w2$:\,for 
every finite $u \subseteq w$ for some $f \in
F$ we have $g\restriction u = f\restriction u \}$. \newline
\noindent  If $f \in {}^w2$, $w \subseteq \text{Ord}$ and 
$\alpha \in \text{ Ord}$
let $f^{[\alpha]}$ be $(f\restriction (w \cap \alpha )) \cup 0_{w \backslash
\alpha}$; let $f^{[\infty]} = f$.\newline
2) Let $\mu = \mu^{<\mu} < \lambda$.  We define a forcing notion
$Q = Q_{\lambda,\mu}$:

\roster
\item "{\bf(a)}"  the members are pairs $p = (w,F) = (w^p,F^p)$, \,
$w \subseteq \lambda$, $|w| < \mu$, and $F$ is a family of $< \mu$ functions
from $w$ to $\{ 0,1 \}$ satisfying
{\roster
\itemitem { $(\alpha)$ }  for every $\alpha \in w$, for some $f \in F$,
$f(\alpha) = 1$
\itemitem { $(\beta)$ }  if $f \in F$ and $\alpha \in w$ then
$f^{[\alpha]} \in F$
\endroster}
\item "{\bf(b)}"  the order: $p \le q$ iff $w^p \subseteq w^q$ and
{\roster
\itemitem { $(\alpha)$ }  $f \in F^q \Rightarrow f\restriction w^p \in c
\ell(F^p)$
\itemitem { $(\beta)$ }  $\forall f \in F^p\exists g \in F^q(f \subseteq g)$.
\endroster}
\endroster
\medskip

\noindent
3)  For $w \subseteq \lambda$, $F \subseteq {}^w2$ let
$BA[w,F] = BA[(w,F)]$
be the Boolean algebra freely generated by $\{ x_\alpha:\alpha \in w \}$
except that: if $u$, \, $v$ are finite subsets of $w$ and for no $f \in F$,
$1_u \cup 0_v \subseteq f$ then $\dsize \bigcap \Sb {\alpha \in u} \endSb
x_\alpha - \bigcup \Sb {\beta \in v} \endSb x_\beta = 0$.
\newline
4)  If $G \subseteq Q_{\lambda,\mu}$ is generic over $V$ then
$BA[G]$ is $\bigcup \Sb {p \in G} \endSb BA[p]$ (see 2.6(2),(3) below).
Here $BA[p] =: BA[w^p,F^p]$.
\enddefinition  
\bigskip

\proclaim{2.6 Claim}  0) For $p \in Q_{\lambda,\mu}$, $BA[p]$ is a Boolean
algebra; also for $f \in F^p$ and ordinal $\alpha \in w^p$ 
(or $\infty$) we have
$f^{[\alpha]} \in F^p$. \newline
1) If $f \in F^p$, $p \in Q_{\lambda,\mu}$ then $f$ induces a homomorphism we
call $f^{\text{hom}}$ from $BA[p]$ to the two members Boolean algebra 
$\{ 0,1 \}$.
In fact for a term $\tau$ in $\{ x_\alpha:\alpha \in w^p \}$, $BA[p] \models
``\tau \ne 0"$ \underbar{iff} for some $f \in F^p$, $f^{\text{hom}}(\tau) = 
1$. \newline
2)  If $p \le q$ then $BA[p]$ is a Boolean subalgebra of $BA[q]$. \newline
3)  Hence $BA[\underset \sim\to G]$ is well defined, $p \Vdash$ ``
$BA[p]$ is a Boolean subalgebra of $BA[\underset \sim\to G]$''. \newline
4)  For $p \in Q_{\lambda,\mu}$, for $\alpha \in w^p,x_\alpha$ 
is a non-zero element which
is not in the subalgebra generated by $\{x_\beta:\beta < \alpha \}$ nor is
there below it a non-zero member of \newline
$\langle x_\beta:\beta < \alpha \rangle_{BA[p]}$.
\endproclaim
\bigskip

\demo{Proof}  Check.  Part (0) should be clear, also part (1).  
Now part (2) follows by 2.5(2)(b) 
and the definition of $BA[p]$; so (3) should become clear.
Lastly, concerning part (4), $x_\alpha$ is a non-zero member of $BA[p]$ by
clause $(\alpha)$ of 2.5(2)(a).  For $\alpha \in w^p$, by 2.5(2)(a)$(\alpha)$
there is $f^1 \in F^p$, $f^1(\alpha) = 1$, and by 2.5(2)(a)$(\beta)$ there is
$f^0 \in F^p$, $f^0(\alpha) = 0$, $f^0 \restriction (w \cap \alpha) = f^1
\restriction (w \cap \alpha)$; together with part (1) this proves the second
phrase of part (4).  As for the third phrase, let $\tau$ be a non-zero
element of the subalgebra generated by $\{ x_\beta:\beta < \alpha \}$, so for
some $f \in F^p$, $f^{\text{hom}}(\tau) = 1$.  By 2.5(2)(a)($\beta$),
letting $f_1 = f^{[\alpha]}$, we have $f_1(\alpha) = 0$ and
$f_1 \in F^p$ and $f_1 \restriction (w \cap \alpha) \subseteq f$.  Hence
$f^{\text{hom}}_1(\tau) = f^{\text{hom}}(\tau) = 1$ and 
$f_1(\alpha) = 0$, hence $f^{\text{hom}}_1(x_\alpha) = 0$.  
This proves $BA[p] \models ``\tau \nleq x_\alpha"$.
\hfill$\square_{2.6}$
\enddemo
\bigskip

\proclaim{2.7 Claim}  Assume $\mu = \mu^{<\mu} < \lambda$.
\roster
\item  $Q_{\lambda,\mu}$ is a
 $\mu$-complete forcing notion of cardinality $\le \lambda^{<\mu}$.
\item  $Q_{\lambda,\mu}$ satisfies the $\mu^+$-c.c.
\endroster
\endproclaim
\bigskip

\demo{Proof of 2.7}  1) The number of elements of $Q_{\lambda,\mu}$ is at 
most

$$
\align
|\{(w,F):&w \subseteq \lambda,|w| < \mu \text{ and } F \text{ is a family of }
< \mu \text{ functions from} \\
  &w \text{ to } \{0,1\} \}| \\
  &\le \sum \Sb w \subseteq \lambda \\ |w| < \mu \endSb
|\{F:F \subseteq {}^w 2, \text{ and } |F| < \mu\}| \\
  &\le \sum \Sb w \subseteq \lambda \\ |w| < \mu \endSb
(2^{|w|})^{< \mu} \le |\{w:w \subseteq \lambda \text{ and }
|w| < \mu\}| \times \mu \\
  &= \lambda^{< \mu} + \mu = \lambda^{< \mu}.
\endalign
$$
\medskip

\noindent
As for the $\mu$-completeness, let $\langle p_\zeta:\zeta < \delta \rangle$ be
an increasing sequence of members of $Q_{\lambda,\mu}$ with $\delta < \mu$.
Let $p_\zeta = (w_\zeta,F_\zeta)$, let $F'_\zeta = c \ell(F_\zeta)$, let
$w = \dsize \bigcup_{\zeta < \delta} w_\zeta$ and let $F' =
\{ f \in {}^w 2:\text{for every } \zeta < \delta \text{ we have }
f \restriction w_\zeta \in F'_\zeta\}$. \newline
Clearly for every $\zeta < \delta$ and $f \in F_\zeta$ there is
$g = g_f \in F'$ extending $F$.  Lastly, let $F = \{g_f:f \in
\dsize \bigcup_{\zeta < \delta} F_\zeta\}$.  Then $p = (w,F) \in Q_{\lambda,
\mu}$ is an upper bound of $\langle p_\zeta:\zeta < \delta \rangle$, as 
required. \newline
2) By the $\Delta$-system argument it suffices to prove that $p^0$,
$p^1$ are compatible when:
\roster
\item "{\bf(a)}"  $otp\left(w^{p^0} \right) = otp\left(w^{p^1} \right)$
and (letting $H = H^{OP}_{w^{p^1},w^{p^0}}$ be the unique
order preserving function from $w^{p^0}$ onto $w^{p^1}$),

\item "{\bf(b)}"  $H$ maps $p^0$ onto $p^1$ i.e.
$$
f \in F^{p^0} \Leftrightarrow (f \circ H^{-1}) \in F^{p^1}
$$

\item "{\bf(c)}"  $\alpha \in w^{p^0} \Rightarrow \alpha \le H(\alpha)$
\item "{\bf(d)}"  for $\alpha \in w^{p^0}$ we have $\alpha \in w^{p^1}$ iff
$\alpha = H(\alpha)$.
\endroster

We define now $q \in Q:w^q = w^{p^0} \cup w^{p^1}$,
$$
F^q = \left\{ \left( f \cup (f \circ H) \right)^{[\beta]}:
f \in F^{p^1},\beta \in w^q \cup \{ \infty \} \right\}.
\qquad\hfill \square_{2.7}
$$
\enddemo
\bigskip

\proclaim{2.8 Claim} Suppose $Q = Q_{\lambda,\mu}$ and
\roster
\item "{\bf(a)}"  $p^\ell \in Q$ for $\ell < m$
\item "{\bf(b)}"  otp$\left( w^{p^\ell} \right) = \text{ otp}(w^{p^0})$, and
$H_{\ell ,k} = H^{OP}_{{w^{p^\ell}},{w^{p^k}}}$ (see the proof of 2.7(2))
\item "{\bf(c)}"  $H_{\ell,k}$ maps $p^k$ onto $p^\ell$,
\item "{\bf(d)}"  for $\alpha \in w^{p^0}$ the sequence $\langle H_{\ell,0}
(\alpha):\ell = 1,\dotsc,m-1 \rangle$ is either strictly increasing 
or constant; and $\{ \alpha,\beta \} \subseteq w^{p^0} \and \ell$, $k < m \and
 H_{\ell,0}(\alpha) = H_{k,0}(\beta)$ implies
$\alpha = \beta$.  Lastly letting $w^\ast = w^{p^0} \cap w^{p^1}$ we have
\newline
$[\ell \ne k \Rightarrow w^{p^\ell} \cap w^{p^k} = w^\ast]$ and $H_{\ell,k}
\restriction w^\ast$ is the identity.
\item "{\bf(e)}"  $\tau(x_1,\ldots ,x_n)$ is a Boolean term, $\alpha^0_i \in
w^{p^0}$ for \newline
$i \in \{1,\dotsc,n\},\alpha^0_1 < \cdots < \alpha^0_n$,
$\alpha^\ell_i = H_{\ell,0}(\alpha^0_i)$.
\item "{\bf(f)}"  In $BA[p^0]$, $\tau \left( x_{\alpha^0_1},\dotsc,
x_{\alpha^0_n} \right)$ is not zero and even not in the subalgebra generated
by $\{ x_\alpha:\alpha \in w^\ast \}$.
\item "{\bf(g)}"  $m-1 > n+1$.
\endroster
\underbar{Then} there is $q \in Q$ such that:
\roster
\item "{\bf($\alpha$)}"  $p^\ell \le q$ for $\ell < m$; and $w^q = \dsize
\bigcup_{\ell <m}w^{p_\ell}$
\item "{\bf($\beta$)}"  $q \Vdash$ ``in $BA[\underset\sim\to G]$,
there is a non-zero Boolean combination $\tau^*$ of \newline
$\left\{ \tau \left(
x_{\alpha^\ell_1},\dotsc,x_{\alpha^\ell_n} \right):1 \le \ell < m 
\right\}$ which is $\le \tau \left(
x_{\alpha^0_1},\dotsc,x_{\alpha^0_n} \right)$".
\endroster
\endproclaim
\bigskip

\demo{Proof}  By assumption (f) (and 2.6(0),(4)) there are $f^\ast_0$, \,
$f^\ast_1 \in c\ell \left( F^{p^0} \right)$ such that:
\roster
\item "{(A)}"  $f^\ast_0 \restriction w^\ast = f^\ast_1 \restriction
w^\ast$
\item "{(B)}"  in the two members Boolean algebra $\{ 0,1 \}$ we have

$$
\tau \left( f^\ast_0(\alpha^0_1),\dotsc,f^\ast_0(\alpha^0_n) \right) = 0
$$

$$
\tau \left( f^\ast_1(\alpha^0_1),\dotsc,f^\ast_1(\alpha^0_n) \right) = 1.
$$
\endroster

\noindent  Now there is $\gamma \in w^{p^0} \cup \{ \infty \}$ such that
$(f^\ast_0)^{[\gamma ]} = f^\ast_0 \and (f^\ast_1)^{[\gamma ]} = f^\ast_1$ \,
(e.g. $\gamma = \infty$).  Choose such $(\gamma,f^\ast_0,f^\ast_1)$ with
$\gamma$ minimal.
\enddemo
\bigskip

\noindent
Let $w^q = \dsize \bigcup_{\ell=0}^{m-1} w^{p^\ell}$. 
We define a function $g \in {}^{(w^q)}2$ as follows:\newline
\demo{\underbar {First Case}}  $g \restriction w^{p^0} = f^\ast_1$ \enddemo
\smallskip

\demo{\underbar {Second Case}}  For odd $\ell \in [1,m)$, \,
$g \restriction w^{p^\ell} = f^\ast_1 \circ H_{0,\ell}$ and \enddemo
\smallskip

\demo{\underbar {Third Case}}  For even $\ell \in [1,m)$, \, (but not
$\ell = 0$!) \, $g^\ast \restriction w^{p^\ell} = f^\ast_0 \circ H_{0,\ell}$.
\enddemo
\medskip

Now $g$ is well defined by clause (A) above.  Let us define $q$:
\medskip
$$
\align
F^q = \biggl\{ \left( \dsize \bigcup _{\ell=0}^{m-1} f \circ H_{0,\ell}
\right) ^{[\alpha]}:&\,\alpha \in w^q \cup \{ \infty \} \text{ and } f \in
F^{p_0} \biggr\}  \\
  &\,\bigcup \left\{ g^{[\alpha]}:\alpha \in w^q \cup \{ \infty \}
\right\}.
\endalign
$$

$$
q = (w^q,F^q).
$$
\enddemo
\medskip
\noindent  Let us check the requirements. \newline
\medskip

\subhead {First Requirement: $q \in Q$} \endsubhead

Clearly $w^q \in [\lambda]^{< \mu}$.  Also $F^q \subseteq {}^{(w^q)}2$, 
$|F^q| < \mu$ so we have to check the conditions $(\alpha)$ and $(\beta)$ of
Definition 2.5(2)(a):
\bigskip

\demo{Condition $(\alpha)$}

If $\alpha \in w^q$ then for some $\ell < m$, \, $\alpha \in w^{p^\ell}$, \,
so as $p^\ell \in Q$ there is $f_\ell \in F^{p^\ell}$ such that $f_\ell
(\alpha)=1$.  Now for some $f_0 \in F^{p^0}$, \, $f_0 = f_\ell \circ
H_{0,\ell}$ so \newline
\noindent $f=: \dsize \bigcup _{k<m} (f_0 \circ H_{0,k}) = \left(
\dsize \bigcup _{k<m} (f_0 \circ H_{0,k}) \right)^{[\infty]}$ belongs to $F^q$
and \newline
\noindent $f(\alpha) = (f_0 \circ H_{0,\ell})(\alpha) = f_\ell(\alpha) = 1$.
\enddemo
\medskip

\demo{Condition $(\beta)$}

As for $\alpha, \beta \in w \cup \{ \infty \}$ and $f \in {}^{(w^q)}2$ we have
$\left( f^{[\alpha]} \right)^{[\beta]} = f^{[Min \{ \alpha,\beta \}]}$
and as $\biggl( \dsize \bigcup_{\ell < m} f_\ell \biggr)^{[\alpha]} =
\dsize \bigcup_{\ell < m} (f_\ell)^{[\alpha]}$, this condition
holds by the way we have defined $F^q$.
\enddemo
\medskip

\subhead {Second Requirement: For $\ell < m, \, p^\ell \le q$} \endsubhead 

By the choice of $q$ clearly $w^{p^\ell} \subseteq w^q$.

\noindent  Also if $f \in F^{p^\ell}$ then $(f \circ H_{\ell,0}) \in F^{p^0}$
and
$$
\bigcup _{k<m} \left( (f \circ H_{\ell,0}) \circ H_{0,k} \right)^{[\infty]}
$$
\medskip

\noindent
belongs to $F^q$ and extends $f$.\newline
Lastly, if $f \in F^q$ we shall prove that $f \restriction
w^{p^\ell} \in c \ell(F^{p^\ell})$ (in fact, $\in F^{p^\ell}$); we have 
two cases:
in the first case $f = \left( \dsize
\bigcup _{\ell<m} (f_0 \circ H_{0,\ell}) \right)^{[\alpha]}$ for some $f_0 \in
F^{p^0}$, let $\beta = \text{ Min } [w^\ell \cup \{ \infty \}
\backslash \alpha]$, \,
so $f^{[\alpha]} \restriction w^\ell = (f_0 \circ H_{0,\ell})^{[\beta]}$,
clearly $f_0 \circ H_{0,\ell} \in F^{p^\ell}$ hence $(f_0 \circ H_{0,\ell})
^{[\beta]} \in F^{p^\ell}$ is as required.
The second case is $f = g^{[\alpha]}$, let \newline
$\beta = \text{ Min } \left[ w^{p^\ell} \cup \{ \infty \}
\backslash \alpha
\right]$, now $f \restriction w^{p^\ell}$ is $f^\ast_0 \circ H_{0,\ell}$ or
$f^\ast_1 \circ H_{0,\ell}$ so $f \restriction w^{p^\ell}$ is $(f^\ast_0 \circ
H_{0,\ell})^{[\beta]}$ or $(f^\ast_1 \circ H_{0,\ell})^{[\beta]}$ hence
belongs to $F^{p^\ell}$.
\bigskip

\subhead {Third Requirement: There is a non-zero Boolean combination of
\newline
\noindent $\left\{ \tau \left(
x_{\alpha^\ell_1},\dotsc,x_{\alpha^\ell_n} \right):\ell = 1,m-1 \right\}$
which is $\le \tau \left( x_{\alpha^0_1},\dotsc,x_{\alpha^0_n} \right)$
in $BA[q]$} \endsubhead
\bigskip

The required Boolean combination will be
$$
\tau^\ast = \bigcap_{\ell=0}^{\left[ \frac {m-2}2 \right]} \tau
\left( x_{\alpha^{2 \ell +1}_1},\dotsc,x_{\alpha^{2 \ell +1}_n} \right) -
\bigcup_{\ell=1}^{\left[ \frac {m-1}2 \right]} \tau
\left( x_{\alpha^{2 \ell}_1},\dotsc,x_{\alpha^{2 \ell}_n} \right).
$$
\medskip
\noindent  So we have to prove the following two assertions.

\bigskip

\noindent  \underbar{First assertion}:  $BA[q] \models ``\tau^\ast \ne 0"$.

Now $g = g^{[\infty]} \in F^q$ satisfies, for each $\ell \in 
\left[ 0, \left[ \frac{m-2}2 \right] \right]$:

$$
g^{\text{hom}} \left( \tau \left( x_{\alpha^{2 \ell +1}_1},\dotsc,
x_{\alpha^{2 \ell+1}_n}
\right) \right) = \left( f^\ast_1 \circ H_{0,2\ell +1} \right)^{\text{hom}} 
\left(
\tau \left( x_{\alpha^{2 \ell +1}_1},\dotsc,x_{\alpha^{2 \ell +1}_n} \right)
\right) = 1;
$$
\medskip

\noindent  also for each $\ell \in \left[ 1,\left[ \frac {m-1}2 \right]
\right]$,

$$
g^{\text{hom}} \left( \tau 
\left( x_{\alpha^{2 \ell}_1},\dotsc,x_{\alpha^{2 \ell}_n}
\right) \right) = \left( f^\ast_0 \circ H_{0,2\ell} \right)^{\text{hom}} 
\left(
\tau \left( x_{\alpha^{2 \ell}_1},\dotsc,x_{\alpha^{2 \ell}_n} \right)
\right) = 0.
$$
\medskip

\noindent
Putting the two together, we get the assertion.

\medskip
\noindent  \underbar{Second assertion}:
$BA[q] \models ``\tau \left( x_{\alpha^0_1},
\dotsc,x_{\alpha^0_n} \right) \ge \tau^\ast "$.

So we have to prove just that:
$$
f \in F^q \Rightarrow f^{\text{hom}} \left( \tau^\ast - \tau \left(
x_{\alpha^0_1},\dotsc,x_{\alpha^0_n} \right) \right) = 0.
$$
\medskip
\demo{\underbar {First Case}}  For some $\alpha \in w^q \cup \{ \infty \}$ and
$f_0 \in F^{p^0}$ we have
$$
f = \left( \bigcup_{\ell < m} f_0 \circ H_{0,\ell} \right)^{[\alpha]}.
$$
\medskip
\noindent  Let $\beta_\ell = Min \left( w^{p^\ell} \cup \{ \infty \}
\backslash \alpha \right)$, and let $\gamma_\ell \in w^{p^0}$ be such
that \newline
\noindent 
$\gamma_\ell = H_{0,\ell}(\beta_\ell)$ or $\gamma_\ell = \beta_\ell = \infty$.
\newline

\noindent  Now by the assumption on $\langle w^{p^\ell}:\ell < m \rangle$
we have $\langle \gamma_\ell:\ell < m \rangle$ is non-increasing.  Let for
$\ell < m$, \, $j_\ell = Min \left\{ j:j = n+1 \, \, \text{ or } j \in \{1,
\dotsc,n\} \text{ and } \alpha^0_j \ge \gamma_\ell \right\}$.
So \newline
\noindent $\langle j_\ell:\ell < 
m \rangle$ is non-increasing and there are $\le n+1$ possible values for each
$j_\ell$.  But by
assumption (g), $m-1>n+1$, so for some $k, 0 < k < k + 1 < m$ and
$j_k = j_{k+1}$.\newline
\noindent  So (as $\alpha^i_1 < \dots < \alpha^i_n$)
$$
\bigwedge_{j=1}^n f\left( x_{\alpha^k_j} \right) = f\left( x_{\alpha^{k+1}_j}
\right),
$$
hence
$$ 
\bigwedge _{j=1}^n f^{\text{hom}} \left( x_{\alpha^k_j} \right) = 
f^{\text{hom}}\left( x_{\alpha^{k+1}_j} \right)
$$
hence
$$
f^{\text{hom}} \left( \tau \left(x_{\alpha^k_1},\dotsc,x_{\alpha^k_n} \right)
\right) = f^{\text{hom}} \left( \tau \left(x_{\alpha^{k+1}_1},\dotsc,
x_{\alpha^{k+1}_n} \right) \right), 
$$
\medskip
\noindent  hence (see Definition of $\tau^\ast) \quad f^{\text{hom}}
(\tau^\ast) = 0$ \,\, hence
$$
f^{\text{hom}} 
\left( \tau^\ast - \tau \left(x_{\alpha^0_1},\dotsc,x_{\alpha^0_n}
\right) \right) = 0,
$$
as required.\enddemo

\medskip
\demo{\underbar {Second Case}}  For some $\alpha \in w^q \cup \{ \infty \},
\quad f = g^{[\alpha]}$.

Let again 
$\beta_\ell = \text{ Min} \left( w^{p^\ell} \cup \{ \infty \} \backslash
\alpha \right)$, $\gamma_\ell = H_{0,\ell}(\beta_\ell)$ (or $\gamma_\ell =
\beta_\ell = \infty$), and

$$
\gamma'_\ell = \cases \gamma_\ell \quad 
\text{ if } \quad \gamma_\ell<\gamma \\
  \infty \quad \text{ if } \quad \gamma_\ell \ge \gamma \endcases
$$
\medskip

\noindent
$j_\ell = \text{ Min}\left\{ j:j = n+1 \text{ or } j \in
\{1,\dotsc,n\} \text{ and } \alpha^0_j \ge \gamma'_\ell \right\}$.  
So \newline
$\langle \gamma_\ell:\ell < m \rangle,\langle \gamma'_\ell:\ell < m
\rangle$ are non-increasing and so is $\langle j_\ell:\ell < m \rangle$.
Here $\gamma$ is the ordinal we chose before defining $q_1$ just after (B) in
the proof.

If for some $k, 0<k<k+1<m$, $j_k=j_{k+1} \le n$, \newline (hence
$\gamma'_{k+1} \le \gamma'_k < \gamma$) then \newline
$f^{\text{hom}}\left( \tau^\ast - \tau \left( x_{\alpha^0_1},
\dotsc,x_{\alpha^0_n}
\right) \right) = 0$ as $f^{\text{hom}}(\tau^\ast) = 0$ 
which holds because \newline
$f^{\text{hom}}\left( \tau \left( x_{\alpha^k_1},
\dotsc,x_{\alpha^k_n} \right) \right) = f^{\text{hom}}\left( \tau \left(
x_{\alpha^{k+1}_1},\dotsc,x_{\alpha^{k+1}_n} \right) \right)$ (the last
equality holds by the choice of $\gamma$; i.e. if inequality holds then
the triple $(\gamma_k,(f^\ast _0)^{[\gamma_k]},(f^\ast_1)^{[\gamma_k]})$
contradicts the choice of $\gamma$ as minimal).
But $j_\ell \, (\ell = 1,m-1)$ is non-increasing hence we can show inductively
on $\ell = 1,...,n+1$ that $j_{m-\ell} \ge \ell$.  So necessarily
$j_1 = n+1$ but as $j_\ell$ is non-increasing clearly  $j_0 = n+1$ hence
$$
\split
f^{\text{hom}}\left( \tau 
\left( x_{\alpha^0_1},\dotsc,x_{\alpha^0_n} \right) \right)
&=g^{[\alpha]} \left( \tau \left( x_{\alpha^0_1},\dotsc,x_{\alpha^0_n} \right)
\right) = g \left( \tau \left( x_{\alpha^0_1},\dotsc,x_{\alpha^0_n} \right)
\right)\\
&=f^\ast_1 \left( \tau \left( x_{\alpha^0_1},\dotsc,x_{\alpha^0_n}
\right) \right) = 1
\endsplit
$$
hence
$$
f^{\text{hom}} 
\left( \tau^\ast - \tau \left( x_{\alpha^0_1},\dotsc,x_{\alpha^0_n}
\right) \right) = 0,
$$
as required.\hfill$\square_{2.8}$
\enddemo
\bigskip

\proclaim{2.9 Theorem}  Suppose $\mu = \mu^{<\mu} < \lambda$,
$Q=Q_{\lambda,\mu}$ and $V \models$ G.C.H. (for simplicity).
\newline
\noindent  Then:
\roster
\item  $Q$ is
$\mu$-complete, $\mu^+$ - c.c. (hence forcing with $Q$ preserves cardinals 
and cofinalities).
\item  $\Vdash_Q ``2^\mu = \left( \lambda^\mu \right)^V"$,
$|Q| = \lambda^{<\mu}$, so
 cardinal arithmetic in $V^Q$ is easily determined.
\item  Let $G \subseteq Q$ be generic over $V$.  \underbar{Then}
$BA[G]$ (see Definition 2.5(4)) is a Boolean algebra of 
cardinality $\lambda$ such that:
{\roster
\itemitem{(a)} if $\theta \le \lambda$ is regular \underbar{then} for 
some subalgebra $B$ of $BA[G]$, $\pi(B) = \theta$
\itemitem{(b)}  if $\theta \le \lambda$ and $\theta > \text{ cf}(\theta) > 
\mu$ \underbar{then} for no $B \subseteq BA[G]$ is $\pi(B) = \theta$
\itemitem{(c)}  $BA[G]$ has $\mu$ non-zero pairwise disjoint elements but no
$\mu^+$ such elements (so its cellularity is $\mu$)
\itemitem{(d)}  if $a \in B^+$ then $BA[G] \restriction a$
satisfies (a), (b), (c) above (also (e))
\itemitem{(e)} if $\theta \le \lambda$ and cf$(\theta) \le \mu$ 
\underbar{then}
for some $B^\prime \subseteq BA[G]$ we have $\pi(B^\prime) = \theta$
\itemitem{(f)}  in $BA[G]$ for every $\alpha < \lambda$, $\{ x_\beta:
\alpha \le \beta < \alpha+\mu^+ \} \subseteq B^+$ is
dense in $\langle \{ x_\beta:\beta < \alpha \} \rangle_{BA[G]}$.
\endroster}
\endroster
\endproclaim
\bigskip

\remark{2.9A Remark}  1) This shows the consistency of a negative answer to
problems 13 + 15 of Monk \cite{M}. \newline
\noindent
2) We could of course make $2^\mu$ bigger by adding the right number of 
Cohen subsets of $\mu$.
\endremark
\bigskip

\demo{Proof}  By Claim 2.7 clearly parts (1), (2) hold.  We are left with
part (3), by 2.6(3) $BA[G]$ is a Boolean algebra, by 2.6(4) it has cardinality
$\lambda$.  As for clause (a), it is exemplified  by $\langle x_\alpha:
\alpha < \theta \rangle_{BA[G]}$ (by 2.6(4)).  Clause (c), the first
statement is easy by the
genericity of $G$ (i.e. as for $p \in Q$, $\alpha \in \lambda \backslash
w^p$ we can find $q$, $p \le q \in Q$, $w^q = w^p \cup \{ \alpha \}$ and
in $BA[q]$, $x_\alpha$ is disjoint to all $y \in J$, for any ideal $J$ of
$BA[p])$.  As for clause (c), second statement, it follows from the
$\triangle$-system argument and the proof of 2.7(2). \newline

Clause (d): 
concerning the generalization of clause (a), let $a \in (BA[G])^+$,
so we can find finite disjoint $u,v \subseteq \lambda$ such that
$0 < \dsize \bigcap_{\alpha \in u} x_\alpha - \dsize \bigcap_{\alpha \in v}
x_\alpha \le a$, choose $\beta = \sup(u \cup v)$, let 

$$
\align
U = \biggl\{ x_\gamma:&\beta < \gamma < \beta + \mu, \text{ and} \\
  &BA[G] \models ``x_\gamma \le \left( \dsize \bigcap_{\alpha \in u}
x_\alpha - \dsize \bigcap_{\alpha \in v} x_\alpha \right)" \biggr\}.
\endalign
$$
\medskip

\noindent
This set is forced to be of cardinality $\mu$ and the subalgebra of
$(BA[G]) \restriction a$ generated by 
$\{x_\gamma \cap a:\gamma \in [\beta,\beta + \mu)\}$ is as required. \newline
Clause (d), the generalization of clause (b) will follow from clause (b).
Clause (d) the generalization of clause (c); the cellularity $\le \mu$ follows
from clause (c), the existence of $\mu$ pairwise disjoint elements follows
from: for every $p \in Q_{\lambda,\mu},\alpha < \lambda$ and $a \in BA[p]$
such that $a \in \langle x_\beta:\beta \in w^p \cap \alpha \rangle_{BA[p]}$
and $\beta \in [\alpha,\lambda) \backslash w^p$ there is $q$ such that
$p \le q \in Q_{\lambda,\mu}$ and $BA[q] \models ``\dsize \bigwedge
_{\gamma \in w^p \backslash \alpha} x_\beta \cap x_\gamma = 0 \and x_\beta \le
a"$.

As for clause (e), (and the generalization in clause (d)), let
$a \in (BA[G])^+$, let $u \subseteq \lambda$ be finite such that
$a \in \langle x_\alpha:\alpha \in u \rangle_{BA[G]}$.  Then we can find
$\langle a_i:i < \mu \rangle$ pairwise disjoint non-zero member of
$\langle x_\alpha:\alpha \in [\sup(u),\sup(u) + \mu)\rangle_{BA[G]}$ which
are below $a$.  Let $\theta = \dsize \sum_{\zeta < \text{ cf}(\theta)}
\theta_\zeta$, each $\theta_\zeta$ regular, let $B_\zeta \subseteq (BA[G])
\restriction a_\zeta$ be a subalgebra with $\pi(B_\zeta) = \theta_\zeta$, and
lastly let $B$ be the subalgebra of $BA[G] \restriction a$ generated by
$\{ a_i:i < \text{ cf}(\theta)\} \cup \dsize \bigcup_{\zeta < \text{ cf}
(\theta)} B_\zeta$; check that $\pi(B) = \theta$. \newline

Clause (f) follows by a density argument.  The real point
(and the only one left) is to prove Clause (b) of part (3).  So suppose
toward contradiciton that $\mu < \text{ cf}(\chi) < \chi \le \lambda$ and 
$p \in Q$ but
$p \Vdash_Q``{\underset\sim\to B}
\subseteq BA{[\underset\sim\to G]}$ is a subalgebra,
$\pi{(\underset\sim\to B)} = \chi"$.

So by Claim 2.2(1)+2.3(1)\newline
$p \Vdash$ ``for arbitrarily large regular $\theta < \chi$, there
is $\bar y = \langle \bar y_\alpha:\alpha < \theta \rangle$\newline

(a sequence of non-zero elements of $\underset\sim\to B$) and
$\theta$-complete proper filter $D$ on $\theta$\newline

(containing the cobounded subsets of $\theta$) such that $\otimes^B_{{\bar y},
D}$ holds (see 2.3(1))''.

\noindent  Let $\kappa = \text{ cf}(\chi)$, so we can find 
regular $\theta_\zeta \in (\text{cf}(\chi),\chi)$, 
(so $\theta_\zeta > \mu)$ increasing with $\zeta$,
$\chi = \dsize \sum_{\zeta < \kappa} \theta_\zeta$, for $i < \kappa$, 
$\biggl( \dsize \sum_{j < i} \theta_j
\biggr)^\kappa < \theta_i$ (remember $V \models$ G.C.H.) and for each
$\zeta < \kappa$, condition $p_\zeta,p \le p_\zeta \in Q$,
and
${\underset\sim {}\to {\bar y}^\zeta} =
\langle {\underset\sim {}\to {\bar y}_\alpha^\zeta}
:\alpha < \theta_\zeta \rangle$, and ($Q$-name of a) proper
 $\theta_\zeta$-complete filter $\underset\sim\to D_\zeta$ on
$\theta$ containing the co-bounded subsets of $\theta$ 
such that
$p_\zeta \Vdash ``\otimes^B_{\underset\sim {}\to {\bar y}^\zeta, 
\underset\sim {}\to D_\zeta}"$
(and without loss of generality
$\Vdash ``\underset\sim\to y^\zeta_\alpha \in (BA[\underset\sim\to G])^+"$).
For each $\zeta < \kappa$ and $\alpha < \theta_\zeta$ there is a
maximal
antichain $\bar p^{\zeta,\alpha} = \langle p_{\zeta,\alpha,\epsilon}:\epsilon
< \mu \rangle$ of members of $Q$ above $p_\zeta$ and terms
$\tau_{\zeta,\alpha,\epsilon} =
\tau^\prime_{\zeta,\alpha,\epsilon}\left( x_{\beta(\zeta,\alpha,
\epsilon,0)},\dotsc,
x_{\beta(\zeta,\alpha,\epsilon,n_{\alpha}(\zeta,\epsilon))} \right)$ 
(i.e. Boolean
terms in $\{ x_\alpha:\alpha < \lambda \}$) such that: \,\, $p_\zeta \le p_
{\zeta,\alpha,\epsilon}$ and $p_{\zeta,\alpha,\epsilon} \Vdash_Q
 ``{\underset \sim\to y^\zeta_\alpha} = \tau_{\zeta,\alpha,\epsilon}"$.
 \newline
Without loss of generality $\{ \beta(\zeta,\alpha,\epsilon,\ell):\ell \le
n_\alpha(\zeta,\epsilon) \} \subseteq w[p_{\zeta,\alpha,\epsilon}]$.
\newline
\noindent
Clearly for each $\zeta < \kappa$, $p_\zeta \Vdash ``\theta_\zeta$ is the
disjoint union of
$\{ \alpha < \theta_\zeta:p_{\zeta,\alpha,\epsilon} \in
{\underset\sim\to G} \}$ for $\epsilon < \mu"$  so for
some $Q$-name ${\underset\sim {}\to \epsilon_\zeta} < \mu$,
we have $p_\zeta \Vdash_Q ``\left\{ \alpha < \theta:p
_{\zeta,\alpha,{\underset\tilde {}\to \epsilon_\zeta}}
\in {\underset\sim {}\to G} \right\}
\ne \emptyset$ mod $\underset\sim \to D_\zeta"$.
So there are $\epsilon_\zeta < \mu$ and $q_\zeta$ satisfying $p_\zeta 
\le q_\zeta \in
Q$, such that $q_\zeta \Vdash ``\epsilon_\zeta$ is as above" and let
$p_{\zeta,\alpha} =: p_{\zeta,\alpha,\epsilon_{\zeta}}$.  So we have
a $Q$-name
$\underset\sim\to A_\zeta$ such that $q_\zeta \Vdash_Q
``\underset\sim\to A_\zeta \subseteq \theta_\zeta$,
$\underset\sim\to A_\zeta \ne \emptyset$ mod $\underset
\sim\to D_\zeta$ and $\alpha \in {\underset\sim\to A_\zeta}
\Leftrightarrow p_{\zeta,\alpha} \in \underset\sim\to G_Q"$.

By possibly replacing $\theta_\zeta,\underset\sim\to A_\zeta$ by
$A^\ast_\zeta \in [\theta_\zeta]^{\theta_\zeta}$, $\underset\sim\to A^\prime
_\zeta = \underset\sim\to A_\zeta \cap A^*_\zeta$ respectively, and
increasing $q_\zeta$, we can assume:\newline
\noindent otp$\left( w^{p_{\alpha,\zeta}} \right) = i^\ast_\zeta(< \mu)$, 
and letting
$w^{p_{\alpha,\zeta}} = \left\{ \beta_{\alpha,\zeta,i}:i < i^\ast_\zeta
\right\}$, (increasing with $i$) and (by Claim 2.4) for some $w^\ast_\zeta
\subseteq i^\ast_\zeta$, $\langle \beta^\ast_{\zeta,i}:i < i^\ast_\zeta
\rangle$ and $\tau^\ast_\zeta$ we have:
$\tau^\prime_{\zeta,\alpha,\epsilon_\zeta} = \tau^\ast_\zeta$,
and for some strictly increasing $\langle j(\zeta,\ell):\ell \le n_\zeta
\rangle$ we have:

$$
\align
q_\zeta \Vdash_Q ``&(a)\,\, \alpha \in \underset\sim\to A_\zeta \and
i \in w^\ast_\zeta \Rightarrow \beta_{\alpha,\zeta,i}
= \beta^\ast_{\zeta,i}\\
                   &(b)\,\, \beta(\zeta,\alpha,\epsilon_\zeta,\ell)
= \beta_{\zeta,\alpha,j(\zeta,\ell)} \text{ and } n_\alpha(\zeta,
\epsilon_\zeta) = n_\zeta\\
                   &(c)\,\, \text{ for every } \beta^\prime_{\zeta,i}
< \beta^\ast_{\zeta,i} (\text{ for } i \in i^\ast_\zeta \backslash w^\ast
_\zeta) \text{ we have }\\
&\quad \biggl\{ \alpha < \theta_\zeta:\alpha \in \underset\sim\to A_\zeta,
[i \in w^*_\zeta \Rightarrow \beta_{\alpha,\zeta,i} = \beta^*_{\zeta,i}]
\text{ and for } i \in i^*_\zeta \backslash w^*_\zeta \\
&\quad \quad \quad \quad \quad \quad\text{ we have } \beta^\prime
_{\zeta,i}< \beta_{\alpha,\zeta,i} < \beta^\ast_{\zeta,i} \biggr\}
\ne \emptyset \text{ mod } \underset\sim\to D_\zeta".
\endalign
$$

\noindent  Also $[i \in i^\ast_\zeta \backslash w^\ast_\zeta \Rightarrow
\biggl( 2 + \dsize \sum_{j<i} \theta_j \biggr)^\kappa < cf(\beta^\ast
_{\zeta,i}) \le \theta_\zeta]$ (remember $\underset\sim\to D
_\zeta$ is a $\theta_\zeta$-complete filter) on $\theta_\zeta$.

As we can replace $\langle \theta_\zeta:\zeta < \kappa \rangle$ by any
subsequence of length $\kappa$, and $\kappa = cf(\kappa) > \mu$ without loss
of generality: $i^\ast_\zeta = i^\ast$, $w^\ast_\zeta = w^\ast,
\text{ otp}(w^{q_\zeta}) = j^\ast$.  Let \newline
$w^{q_\zeta} = \{ \beta^\ast_{\zeta,i}:i^\ast \le i < j^\ast\}$.
Now we apply 2.4 to $\kappa$, $D_\kappa$ (filter of closed unbounded
sets) and $\langle \langle \beta^\ast_{\zeta,i}:i < j^\ast \rangle:\zeta
< \kappa \rangle$ and get $\langle \beta^\otimes_i:i < j^\ast \rangle$ and
$w^\otimes \subseteq j^\ast$.  Without loss of generality the $q_\zeta$
are pairwise isomorphic.  Note \newline
$[i \in i^\ast \backslash w^\ast \and \zeta
< \xi < \kappa \Rightarrow \beta^\ast_{\zeta,i} \ne \beta^\ast_{\xi,i}]$
(as $cf(\beta^\ast_{\zeta,i}) \le \theta_\zeta$,
$cf(\beta^\ast_{\xi,i}) > \theta_\zeta)$. \newline
Hence $w^\otimes \cap i^\ast \subseteq w^\ast$.  For every $\zeta < \kappa$
and $i \in i^\ast \backslash w^\ast$, let $\beta^-_{\zeta,i} <
\beta^\ast_{\zeta,i}$ be such that the interval
$[\beta^-_{\zeta,i},\beta^\ast_{\zeta,i})$ is
disjoint to $\{ \beta^\ast_{\xi,j}:\xi < \kappa,j < i^\ast \} \cup \{ \beta
^\otimes_i:i < i^\ast \}$, and as we can omit an initial segment of
$\langle \theta_i:i < \kappa \rangle$ without loss of generality
$[\beta^\ast_{\zeta,i},\beta^\otimes_i)$ is disjoint to
$\{\beta^\otimes_j:j < i^\ast \}$.  Choose for each $\zeta < \kappa$,
$\alpha_\zeta \in A^\ast_\zeta$, such that
$\biggl[ i \in i^\ast \backslash w^\ast \Rightarrow \beta_{\zeta,\alpha_\zeta,
i} \in [\beta^-_{\zeta,i} + \mu,\beta^\ast_{\zeta,i}) \biggr]$.
Let \newline
\noindent ${\underset\sim {}\to Y}$ = the Boolean subalgebra generated by 
$\left\{ \tau_{\zeta,\alpha_\zeta}:q_\zeta \in \underset\sim\to G
\text{ and } p_{\zeta,\alpha_\zeta} \in \underset\sim\to G \right\}$.
This set has cardinality $\le \kappa$, and we shall prove
$$
q_0 \Vdash_Q ``\underset\sim\to Y \backslash \{ 0 \}
\text{ is dense in } \left\{ \tau_{0,\beta}:\beta \in
\underset\sim\to A_0 \right\}".  \tag$*$
$$
This contradicts the choice $q_0 \Vdash ``(*)_{{\bar y^0},D_0}
^{{\underset\tilde {}\to B}}$, $\underset\sim\to A \ne \emptyset$ mod
$\underset\sim\to D_0"$. \newline
\noindent  To prove $(*)$ assume $q_0 \le r_0 \in Q$, we can
choose $\zeta^\ast < \kappa$ and $r^+_\zeta$, for $\zeta \in
[\zeta^\ast,\kappa)$ such that $r_0 \le r^+_0$, $q_\zeta \le r^+_\zeta$,
$p_{\zeta,\alpha_\zeta} \le r^+_\zeta$ and
$\langle r_\zeta:\zeta = 0 \text{ or }
\zeta \in [\zeta^\ast,\kappa) \rangle$ is as in 2.8; apply 2.8 and
get a contradiction. \hfill$\square_{2.9}$
\enddemo
\bigskip

\noindent
A theorem complementary to 2.9 is:
\proclaim{2.10 Theorem}  If $\pi(B) > \lambda$ and
\roster
\item "{(A)}"  $\text{cf}(\lambda) = \aleph_0$ or
\item "{(B)}"  $\dsize \bigwedge_{\mu < \lambda}\mu^{< \text{ cf}(\lambda)} 
< \lambda$ or
\item "{(C)}"  $\dsize \bigwedge_{\mu < \lambda} 2^\mu < \pi(B)$ (or just)
$\dsize \bigwedge_{\mu < \lambda} 2^\mu < |B| \and \lambda \le \pi(B)$. 
\endroster
\medskip

\noindent
\underbar{Then} $B$ has a subalgebra
$B^\prime$ such that $\lambda = \pi(B^\prime) = |B^\prime|$.
\endproclaim
\bigskip

\remark{Remark}  The conclusion of 2.10 implies that
$\lambda \in \pi_{Ss}(B)$.
\endremark
\bigskip

\demo{Proof}  Case (C) is easier so we ignore it.  
By 2.1 without loss of generality \newline
$\pi(B) = \lambda^+ = |B|$.
We try to choose by induction on $\alpha < \lambda$, $a_\alpha$ such that:
\roster
\item"{(a)}"  $a_\alpha \in B^+$
\item"{(b)}"  for $\beta < \alpha$ we have $B \models ``a_\alpha \cap a_\beta
 = 0"$
\item"{(c)}"  $\pi(B \restriction a_\alpha) < \lambda^+$.
\endroster

Let $a_\alpha$ be defined iff $\alpha < \alpha^\ast$.
\enddemo
\bigskip

\demo{Case a}  $\alpha^\ast \ge \lambda$.\newline
Let $B^\prime$ be the
subalgebra generated by $\{ a_\alpha:\alpha < \lambda \}$ clearly $|B^\prime|
= \pi(B^\prime) = \lambda$.
\enddemo
\bigskip 

\demo{Case b}  Not Case a but $\dsize \sum_{\alpha < \alpha^*}
\pi(B \restriction a_\alpha) \ge \lambda$.\newline

So we can find distinct $\alpha_\zeta < \alpha^\ast$ for $\zeta < 
\text{ cf}(\lambda)$ such that $\dsize \sum_\zeta \pi
\left( B \restriction a_{\alpha_\zeta} \right)
\ge \lambda$.
We can find regular $\theta_\zeta \le \pi \left( B \restriction
a_{\alpha_\zeta}
\right)$, such that $\sup_{\zeta < \text{ cf}(\lambda)} \theta_\zeta = 
\lambda$ then find \newline
$B_\zeta \subseteq B \restriction a_{\alpha_\zeta}$,
such that $|B_\zeta| = \theta_\zeta$ and $\pi(B_\zeta) = \theta_\zeta$ 
(by 2.1).  Let
$B^\prime$ be the subalgebra of $B$ generated by $\dsize 
\bigcup_{\zeta < \text{ cf}(\lambda)}B_\zeta \cup \{ a_{\alpha_\zeta}:
\zeta < \text{ cf}(\lambda)\}$.  \newline
Clearly $|B^\prime| = \pi(B^\prime) = \lambda$.
\enddemo
\bigskip

\demo{Case c}  Cases a, b fail.

Let $I = \{a \in B:\dsize \bigwedge_{\alpha
<\alpha^\ast} a \cap a_\alpha = 0 \}$, so $I$ is an ideal of $B$ and
$a \in I \Rightarrow \pi(B \restriction a) \ge \lambda$.
Also $I \ne \{ 0 \}$ (as if $I = \{0\}$ then 
$\pi(B) \le \dsize \sum_{\alpha < \alpha^\ast}
\pi(B \restriction a_\alpha) < \lambda$).  So easily without
loss of generality:

$$
\text{ if } a \in I \cap B^+ \text{ then } \pi (B \restriction a) > \lambda.
\tag"{$(*)$}"
$$

$$
\text{ if } a \in I \cap B^+ \text{ then } B \restriction a
\text{ is an atomless Boolean algebra}
\tag"{$(**)$}"
$$
\medskip

Now without loss of generality $B$ satisfies (cf$(\lambda)$)-c.c. (otherwise
act as in Case b), so we have finished if Case (A) of the hypothesis
holds. \newline
\noindent So case (B) of the hypothesis holds, hence we can
use Lemma 4.9, p.88 of \cite{Sh:92} and
find a free subalgebra $B^\prime$ of $B$ of
cardinality $(\lambda^{< \text{ cf}(\lambda)})^+$ 
hence of cardinality $\lambda$,
this $B^\prime$ is as required.\hfill$\square_{2.10}$
\enddemo
\newpage

\head {\S3 On $\pi$ and $\pi \chi$ of Products of Boolean Algebras}
\endhead
\bigskip

\proclaim{3.1 Theorem}  Suppose
\medskip
\roster
\item "{$\bigotimes$}"   $\aleph_0 < \kappa = \text{ cf}(\chi) < \chi
< \lambda = \text{ cf}(\lambda) < (\chi^\kappa)$ and
$\dsize \bigwedge_{\theta < \chi} \theta^\kappa< \chi$.
\endroster
\medskip

\noindent
\underbar{Then}
there are Boolean Algebras $B_i$ (for $i < \kappa)$ such that (on $\pi(F,B)$,
$\pi \chi(B)$ see below)
\roster
\item "$(*)$(a)"  $\pi\chi(B_i) < \chi = \dsize \sum_{j < \kappa}
\pi\chi(B_j)$
\item "(b)"  for any uniform ultrafilter $D$ on $\kappa$
$$
\lambda = \pi\chi \left( \dsize \prod_{i<\kappa}B_i/D \right)
$$
\item "(c)"  if $D$ is a regular ultrafilter on $\kappa$ then $\dsize
\prod_{i<\kappa}(\pi\chi(B_i))/D = \chi^\kappa$.
\endroster
\endproclaim
\bigskip

\definition{3.2 Definition}  1) For a Boolean algebra $B$ and ultrafilter $F$
of $B$, let
$$
\pi (F,B) = \text{ Min } \left\{ |X|:X \subseteq B^+
\text{ and } (\forall y \in F)(\exists x \in X)[x \le y] \right\}.
$$
We say $X$ is dense in $F$ (though possibly $x \notin F$).\newline
\noindent 2) For a Boolean algebra $B$,
$$
\pi \chi(B) = \text{ sup } \{ \pi(F,B):F \text{ an ultrafilter of } B \}
$$
$$
\pi \chi^+ (B) = \cup \{ \pi(F,B)^+:F \text{ an ultrafilter of } B \}.
$$
\enddefinition
\bigskip

\remark{3.3 Remark}  1) If $\kappa = \aleph_0$ the theorem holds almost
always and probably always, but we omit it to simplify the statement.
(The theorem holds for $\kappa = \aleph_0$, e.g. if
$\chi < \lambda = \text{ cf}(\lambda) <$ (first fix point $> \chi$),
more generally if
\medskip
\roster
\item "{$\bigotimes'$}"  $\kappa = \text{ cf}(\chi) < \chi < \lambda =
\text{ cf}(\lambda) < pp^+_{J^{bd}_\kappa}(\chi) \text{ and }
2^\kappa < \chi$
\endroster
\medskip

\noindent
see \cite{Sh:g}.  The point is that \cite{Sh:355},5.4 deals with 
uncountable cofinalities).
\endremark
\bigskip

\demo{3.4 Proof of Theorem 3.1}  Let for a linear order ${\Cal I}$, 
$BA[{\Cal I}]$
be the Boolean algebra of subsets of ${\Cal I}$ generated by the close-
open intervals $[a,b) = \{ x \in {\Cal I}:a \le x < b \}$ where we allow
$a \in \{ - \infty
\} \cup {\Cal I}$, $b \in {\Cal I} \cup \{ \infty \}$, (and $a \le b)$.  Now
clearly
\smallskip
\roster
\item "${(*)_1}$" if $F$ is an ultrafilter on ${\Cal I}$,
\underbar{then} there is a
Dedekind cut $({\Cal I}^d,{\Cal I}^u)$ of ${\Cal I}$
(i.e. ${\Cal I}^d \cap {\Cal I}^u = \emptyset$,
${\Cal I}^d \cup {\Cal I}^u = {\Cal I} \text{ and } (\forall x_0 \in
{\Cal I}^d)(\forall x_1 \in {\Cal I}^u)[x_0 < x_1])$ such that for $x
\in BA[{\Cal I}]$, $x \in F$ \underbar{iff} for some $a_0 \in {\Cal
I}^d$, $a_1 \in {\Cal I}^u$ we have
$[a_0,a_1) \le x$.
\item "${(*)_2}$"  if ${\Cal I}, F,
({\Cal I}^d,{\Cal I}^u)$ are as above then

$$
\pi(F,BA({\Cal I})) \text{ is }: \cases \text{Max }
\{ \text{cf}({\Cal I}^d),\text{cf}(({\Cal I}^u)^*)\} &\text{ \underbar{if} }
\quad \text{cf}({\Cal I}^d),\text{cf}(({\Cal I}^u)^*) \ge \aleph_0 \\
  \text{cf}({\Cal I}^d) &\text{ \underbar{if} } \quad \text{cf}
(({\Cal I}^u)^*) = 1 \\
  \text{cf}(({\Cal I}^u)^*) &\text{ \underbar{if} } \quad \text{cf}
({\Cal I}^d) = 1 \\
   1 &\text{ \underbar{if} } \quad \text{cf}({\Cal I}^d) = \text{ cf}
(({\Cal I}^u)^\ast)=1 \endcases
$$
\endroster
\medskip

\noindent
note also

$$
\pi \chi(BA({\Cal I})) = \text{ Max} \{ \text{cf}({\Cal I}^d),
\text{cf}(({\Cal I}^u)^*):({\Cal I}^d,{\Cal I}^u) \text{ a Dedekind cut
of } {\Cal I} \}.
$$
\medskip

\noindent
Now by the assumption $\bigotimes$ (and \cite{Sh:g},II,5.4 +
VIII,\S1]), we can find a (strictly) increasing
sequence $\langle \lambda_i:i < \kappa \rangle$ of regular cardinals,
$\kappa < \lambda_i < \chi$, $\chi = 
\dsize \sum_{i<\kappa} \lambda_i$ such that
$\dsize \prod_{i<\kappa} \lambda_i/J^{bd}_\kappa$ has true
cofinality $\lambda$  (where $J^{bd}_\kappa$ is the ideal of
bounded subsets of $\kappa$).

Let ${\Bbb Q}$ be the rational order and ${\Cal I}_i$ be $\lambda_i
\times {\Bbb Q}$ (i.e. the set of elements is\newline
$\{ (\alpha,q):\alpha < \lambda_i,q \in {\Bbb Q} \}$,
the order is lexicographical).  Let $B_i = BA[{\Cal I}_i]$,
so by $(*)_1,(*)_2$ we know that $\pi \chi(B_i) = \lambda _i$.
Moreover, if $F$ is an ultrafilter of $B_i$, then $\pi(F,B) = \aleph_0$ 
except
when $F$ is the ultrafilter $F_i$ generated by $\{ x^i_\alpha = [(\alpha,0),
\infty):\alpha < \lambda _i \}$.  Let $x^i_{\alpha,q} =: [(\alpha,q),\infty)$.
Let $D$ be a uniform ultrafilter on $\kappa$, so $\dsize \prod_{i<\kappa}
\lambda_i/D$ has cofinality $\lambda$.  Also if $D$ is regular, then (see
\cite{CK}) we know  $\chi^\kappa = \dsize \prod
_{i<\kappa} \lambda_i/D = \dsize \prod_{i<\kappa}(\pi \chi B_i)/D$.
So parts (a) and (c) of $(*)$ of Theorem 3.1 are satisfied.  
To prove part (b) of $(*)$, let $D$ be a
uniform ultrafilter on $\kappa$ and let $B =: \dsize
\prod _{i<\kappa}B_i/D$.  Let $F$ be such that $(B,F) = \dsize \prod _{i<
\kappa}(B_i,F_i)/D$.
Clearly $F$ is an ultrafilter of $B$, it is generated by
$X = \dsize \prod _{i<\kappa} X_i/D$ where $X_i = \{ x^i_{\alpha,q}:\alpha <
\lambda,q \in {\Bbb Q} \} \subseteq B_i$, which is linearly ordered in $B$,
 and this
linear order has the same cofinality as $\dsize \prod _{i<\kappa} \lambda
_i/D$, which has cofinality $\lambda$.  So $\pi \chi(F,B) = \lambda$ hence
$\pi \chi(B) \ge \lambda$.
Let $F^\prime$ be an ultrafilter of $B$, $F^\prime \ne F$. \newline
Let $X_d := \{ x \in X:x \in F^\prime \}$, and \newline
$X_u := \{ x \in X:x \notin F^\prime$
(i.e. $1_B - x \in F^\prime) \}$.  Clearly $(X_d,X_u)$ is a Dedekind
cut of $X$ (which is linearly ordered: as a subset of $B$, or as $\dsize \prod
_{i<\kappa} X_i/D$ where $X_i \subseteq B_i$ inherit the order from $B_i$,
so
$x^i_{\alpha,a} < x^i_{\beta,b} \Leftrightarrow (\alpha,a) < (\beta,b)$.)
If $X_d = X$ easily
$F^\prime = F$ contradiction, so $X_u \ne \emptyset$.

We shall prove now that $\pi \chi(F,B) \le 2^\kappa$.
If not, we shall choose by
induction on $\zeta < \left( 2^\kappa \right)^+$ a set $Y_\zeta$,
subsets $Z^i_\zeta$  of $\lambda_i$ for $i < \kappa$, increasing continuous
in $\zeta$ and $y_\zeta$ such that:
$$
|Z^i_\zeta| \le 2^\kappa
$$

$$
\xi < \zeta \Rightarrow Z^i_\xi \subseteq Z^i_\zeta
$$

$$
Y_\zeta = \dsize \prod _{i<\kappa}(Z^i_\zeta \times {\Bbb Q})/D
 \backslash \{ 0 \}, \text{ so } |Y_\zeta| \le 2^\kappa
$$

$$
y_\zeta \in F^\prime
$$

$$
y_\zeta \in Y_{\zeta + 1}
$$

$$
(\forall y \in Y_\zeta)(y > 0 \Rightarrow \neg y \le y_\zeta).
$$
\smallskip
\noindent There is no problem in doing this for $i = 0$, let
$Z^i_\zeta = \{0\}$, for $i$ limit let $Z^i_\zeta = \dsize
\bigcup_{\varepsilon < \zeta} Z^i_\varepsilon$.  Now having defined 
$\langle Z^i_\zeta:i < \kappa \rangle$ (hence $Y_\zeta$) let

$$
y_\zeta = \langle y^i_\zeta:i < \kappa \rangle /D,
$$

$$
y^i_\zeta = \dsize \bigcup _{\ell < n_{i,\zeta}}
\left[ x^i_{\alpha_{i,\zeta,2 \ell},
q_{i,\zeta,2 \ell}},x^i_{\alpha_{i,\zeta,2 \ell +1},q_{i,\zeta,2 \ell
+1}} \right)
$$
\medskip

\noindent
where $\langle (\alpha_{i,\zeta,\ell},q_{i,\zeta,\ell}):\ell < 
2n_{i,\zeta} \rangle$ is a strictly increasing sequence of members of 
$Z^i_\zeta \cup
\{ - \infty, + \infty \}$ (we write $- \infty = (- \infty,0),+ \infty =
(\infty,0)$.)  Let \newline
$Z^i_{\zeta + 1} = Z^i_\zeta \cup \{ \alpha_{i,\zeta,
\ell}:\ell < 2n_{i,\zeta} \}$.

For some unbounded ${\Cal U} \subseteq \left( 2^\kappa \right)^+$;
\roster
\item "{(a)}"  $q_{i,\zeta,\ell} = q_{i,\ell}$ and $n_{i,\zeta} = n_i$,
\endroster
\medskip

\noindent
apply 2.4, we get an easy contradiction.\hfill$\square_{3.4}$
\enddemo
\bigskip

\remark{3.4A Remark}  We can similarly analyze (when $B_i = BA[{\Cal I}_i]$)
\newline

$$
\left\{ \pi(F,\dsize \prod_{i< \kappa}B_i/D):F
\text{ an ultrafilter
of } \dsize \prod_{i< \kappa} B_i/D \right\}
 \backslash (2^\kappa)^+ =
$$

$$
\align
\biggl\{ \lambda:&\lambda^\prime_i = \text{ cf}(\lambda^\prime_i) > 2^\kappa
 \text{ and in } {\Cal I}_i \text{ there is a Dedekind cut } \\
  &(X_d,X_u) \text{ such that } (\text{cf}(X_d),\text{cf}(X^*_u)) 
= (\lambda^d_i,\lambda^u_i)
\text{ such that } \\
  &\lambda = \text{ Min } [ \{ \text{cf}(\Pi \lambda^u_i/D),
 \text{cf}(\Pi \lambda^d_i/D) \} \backslash \{ 1 \} ] \biggr\}.
\endalign
$$
\endremark
\bigskip

Note in comparison that by Koppelberg Shelah \cite{KpSh:415},Th.1.1
\proclaim{3.5 Theorem}  Assume $D$ is an ultrafilter on $\kappa$, for $i <
\kappa$, $A_i$ is a Boolean Algebra, $\lambda_i = \pi(A_i)$. Assume the
Strong Hypothesis - \cite{Sh:420},6.2, 
i.e. $pp(\mu) = \mu^+$ for all singulars or
just SCH.  If $2^\kappa < \lambda_i$ (or just $\{ i:2
^\kappa < \lambda_i \} \in D)$ then \newline
\noindent $\pi(\dsize \prod_{i< \kappa} A_i/D)
= \dsize \prod_{i< \kappa} \lambda_i/D$.
\endproclaim
\bigskip

\proclaim{3.6 Claim}  Assume that for $i < \kappa$, $A_i$ is an infinite
Boolean
algebra, $D$ is a non-principal ultrafilter on $\kappa$ and
$A =: \dsize \prod_{i<\kappa} A_i/D$.  If  $n_i < \omega$ for $i < \kappa$ and
$\mu =: \dsize \prod_{i<\kappa} n_i/D$ and $\mu$ is a regular cardinal
then $\pi\chi(A) \ge \mu$. 
\endproclaim
\bigskip

\demo{Proof}  Let $\chi$ be a large enough regular cardinal (i.e. such that
$\kappa,D,A_i,A$ belong to $H(\chi)$).  Let ${\frak C}_i = (H(\chi),\in,
<^\ast)$ and ${\frak C} = \dsize \prod_i{\frak C}_i/D$, so $A$ is a member
of ${\frak C}$.

Clearly $\omega^\ast =: \langle \omega:i < \kappa \rangle/D$ is considered
by ${\frak C}$ a limit ordinal and, from the outside, has a cofinality, which
we call $\lambda$.  Without loss of generality $i < \kappa \Rightarrow
n_i > 2$.
\enddemo
\bigskip

\subhead{The proof is divided to two cases} \endsubhead
\demo{Case 1}  There are no $\mu_0 < \mu$ and $n^0_i < n_i$ such that
$\aleph_0 \le \mu_0 = \left| \dsize \prod_{i<\kappa} n^0_i/D \right|$.

\noindent  We can find $n^\ast_i$ such that:
\smallskip
\roster
\item  $\mu = \dsize \prod_{i<\kappa} n^\ast_i/D$
\item  $\mu = \dsize \prod_{i<\kappa} 2^{(n^\ast_i)^{(n^*_i)}}/D$
\endroster

\noindent  Let for $i < \kappa$, $\langle a^i_k:k < 2^{(n^\ast_i)^{(n^*_i)}}
 \rangle$
be pairwise disjoint non-zero members of $A_i$ with union $1_{A_i}$.  Let
$P^i$ be the Boolean subalgebra generated by $\{ a^i_k:k <
 2^{(n^\ast_i)^{(n^*_i)}} \}$.
Let $R^i =: \{ a^i_k:k < 2^{(n^\ast_i)^{(n^*_i)}} \}$.  Let for
$k < n^\ast_i$, $Q^i_k \subseteq P^i$ be such that: $Q^i_k$ is a set of
$n^\ast_i$ pairwise disjoint non-zero elements of $P^i$ such that if
$\langle b_k:
k < n^\ast_i \rangle \in \dsize \prod_{k<n^\ast_i} Q^i_k$ then
$\cap \{ b_k:k < n^\ast_i \}$ is not zero.  \newline
Let $F^i(x) =: \cup
\{ a^i_k:x \cap a^i_k \ne 0_{A_i} \text{ and }
\ell < k \Rightarrow x \cap a^i_\ell = 0_{A_i}\}$ so the union is on at most
one element and $F^i(x) = 0_{A_i} \Leftrightarrow x = 0_{A_i}$. \newline
\noindent  Let $(A,P,Q,R,F,n^\ast) =: \dsize \prod_{i < \kappa} 
(A_i,P^i,Q^i,R^i,F^i,n^\ast_i)$.  (We consider $Q$ as a two place 
relation).  Note that
\smallskip
\roster
\item "$(*)_1$"  $P$ is a Boolean subalgebra of $A$
\endroster
\smallskip

\roster
\item "$(*)_2$"  If $D$ is a subset of $P^+$ then its density in $A$ is
equal to its density in $P$.
\endroster

\noindent [Why?  If $Y \subseteq A^+$ is dense in $D$, then $\{ F(c):c \in
Y \}$ is a subset of $P^+$ dense in $D$ of cardinality $\le |Y|$, for the
other direction use the same set].

\noindent  Now let us enumerate the members of $n^\ast$ as $\{ k_\alpha:
\alpha < \mu \}$ (no repetitions) we also list the members of
$P^+$ as $\{ c_\alpha:\alpha < \mu \}$.  Now we choose
by induction on $\alpha < \mu$ a member $b_\alpha$ of
$Q_{k_\alpha}$ such that it contains (in $A$) no one among
$\{ c_\beta:\beta < \alpha \}$.  As each $c_\beta$ can ``object'' to at most
one $b \in Q_{k_\alpha}$ (as the candidates are pairwise disjoint) and
$Q_{k_\alpha}$ has cardinality $\mu > |\alpha|$ we can do this. Also by
the choice of the $Q^i$'s there is a filter of $P$ to which $b_\alpha$
belongs for every $\alpha < \mu$, so we are done.
\enddemo 
\bigskip

\demo{Case 2 \footnote{in this case the regularity of $\mu$ is not used}}
There is $\mu_0 < \mu$ and $n^0_i < n_i$ such that
$\aleph_0 \le \mu_0 = \dsize \prod_{i < \kappa} n^0_i/D$.
\newline
We can define $X_i, Y_i$ such that: $X_i$ is the family of those
subsets of $Y_i$ with exactly $n^0_i$ elements and $|Y_i| = n_i \times
n^0_i+1$ and e.g. $Y_i$ is a set of natural numbers; note that $|X_i| > n_i$.
Let $X =: \langle X_i:i < \kappa \rangle /D$, 
$Y =: \langle Y_i:i < \kappa \rangle/D$, for $y \in Y$ (in ${\frak C}$'s
sense) let $S_y =: \{ x \in X:y \in x \}$.   Let $n^1_i =: |X_i|$,
note $\left| \dsize \prod_{i<\kappa} n^1_i/D \right| \ge \mu$;
let  $\langle a^i_k:k < n^1_i \rangle$ be a partition of $1_{A_i}$ to non-
zero members of $A_i$ and $h_i$ be a one to one
function from $X_i$ onto $R_i =: \{ a^i_k:k < n^1_i \}$ and for
$y \in Y_i$ let $b^i_y =: \bigcup \{ h_i(x):x \in S_y \} \in A_i$; we define
$h,n^1,\langle b_y:y \in Y \rangle \in {\frak C}$ naturally, let
$P^\ast_i$ be the subalgebra of $A_i$ generated by $\{ a^i_k:k < n^1_i \}$
and $P^* = \dsize \prod_{i < \kappa} P_i/D$ as in the other case.  
By a cardinality argument if $k < \omega$,
$[n^0_i > k \and y_0,\dotsc,y_{k-1} \in Y \Rightarrow
A_i \models ``b^i_{y_0} \cap
\dots \cap b^i_{y_{k-1}} \ne 0_{A_i}"]$ hence $\{ b_y:y \in Y \} \subseteq
P^*$ generates a filter of $P^*$. 
Let $D$ be an ultrafilter of $P^\ast$ containing
$b_y$ for $y \in Y$.  If $Z \subseteq A \setminus \{ 0 \}$ exemplifies the
density of $D$ in $A^\ast$ and is of cardinality $\mu_2 < \mu$, as in case 1
without loss of generality
$Z = \{ a_f:f \in F \} \subseteq P^\ast$ where
$F \subseteq \dsize \prod_{i<\kappa} n^1_i$,
$|F| < \mu$, $a_f =: \langle a^i_{f(i)}:i < \kappa \rangle/D$.  Let
$n = \langle n_i:i < \kappa \rangle/D$, $n^0 = \langle n^0_i:i < \kappa
\rangle/D$.
\newline
Each $f \in F$, $h^{-1}(f)$ is from $X$, so is a subset of $Y$
of cardinality $n_0$ from inside (``considered" by ${\frak C}$ to be so)
so $\mu_0$ from the outside; from the inside has cardinality $n$ and from
the outside $n$ has cardinality
$\mu$ so there is a member of $X$ disjoint to all of the $h^{-1}(f)$,
contradiction to the
density.  So $D$ has in $A$ density $\mu$, hence for every ultrafilter $F$
of $A$ extending $D$, $\pi(F,A) \ge \mu$.  \newline
Hence $\pi\chi(A) \ge \mu$ as required. \hfill$\square_{3.6}$
\enddemo
\bigskip

\remark{Remark}  1) If we ignore regularity, case (1) suffices as for every
$\bar n/D \in \omega^\kappa/D$, \newline
$\mu = \Pi \bar n/D \ge \aleph_0$ we can find
${\bar n^\ell}/D \in \omega^\kappa/D$ such that $\omega^\kappa/D \models
``2^{\bar n^{\ell +1}/D} < {\bar n^\ell}/D"$, so $\langle |{\bar n^\ell}/D|:
\ell < \omega \rangle$ is eventually constant. \newline
2) If each $A_i$ is of cardinality $\aleph_0$, $\mu = \aleph^\kappa_0/D$
is regular the proof above gives $\pi\chi \left( \dsize \prod_{i<\kappa}
A_i/D \right) = \mu$ (if 3.6 does not apply, $\omega^\kappa/D$ is $\mu$-like,
so we can apply case 2 with $|X_i| = |Y_i| = \aleph_0$).
\endremark

\newpage

\vglue -.30truein

\noindent {ADDED IN PROOF}
\bigskip

\noindent We can add to claim 2.3:
\medskip

\noindent {\bf CLAIM 2.3(3):}  Assume that $B$ is a Boolean algebra
and $\sigma<\theta=cf(\theta)$ and $(*)^B_{\sigma} [Y] $ and for no $\tau
\in (\sigma,\theta)$ and $Y' $ do we have $(*)^B_{\tau} [Y']$.

Then in the conclusion of 2.3(1) we can add: $D$ is normal (hence in 2.2 we
get that for arbitrarily large $\theta<\mu$, there are a normal filter $D $
on $\theta$ and $\langle y_i:i<\theta\rangle$ as in 2.3(1) ).
\medskip

\noindent {\bf Remark:}  If $\theta = \tau^+$ then 2.3(2) gives the
conclusion.
\medskip

\noindent {\bf Proof:}
Let $Y= \{y_i:i<\theta\}$, we choose by induction of $n<\omega $ a club
$E_n$ of $\theta$ and sequence $\overline y^n=\langle y^n_i:i<\theta\rangle
$ of non zero members of $B$ such that:

\item {(a)} letting $Y_n=\{y^n_i:i<\theta\}$, we have
$Y=Y_0, Y_n \subseteq Y_{n+1}$ and $E_{n+1} \subseteq E_n$

\item {(b)} if  $\delta\in E_{n+1}$ and
if $\delta <\alpha < \min(E_n \backslash (\delta+1))$ then for some
$\beta<\delta, y^{n+1}_\beta \le_B y^n_\alpha$.

Let for $n=0$, $\overline y^0$ list $Y$ and $E_0=\{\delta<\theta: \delta\;\;
\hbox{a limit ordinal divisible by } \;\;\sigma\}$.

For $n=m+1$, for each $\delta\in E_n$, let $\gamma_\delta =
\min(E_n \backslash (\delta+1))$ and let $Z^n_\delta$ be a subset of $B^+$
of cardinality $\le\sigma$ dense in $\{ y^m_i: \delta \le i <
\gamma_\delta\}$ which exist by the proof of 2.2.  Let $Z^n_\delta = \{
z^n_{\delta + i} : i < \sigma \}$.  (no double use of the same index).

Also for each $\zeta<\sigma$ let $D^n_\zeta$ be the normal filter on
$\theta$ generated by the subsets of $\theta$ of the form $\{ i< \theta: z
\not\le y^n_{i + \zeta}\}$ for $z \in B^+$ ; by our assumption toward
contradiction there are $z^n_{\zeta_\epsilon} \in B^+$ for $\epsilon<\theta$
and club $E^n_\zeta$ of $\theta$ such that if $\delta \in E^n_\zeta$ and
$\zeta<\sigma$ then for some $\epsilon<\delta$ we have
$z_{\zeta,\epsilon}\le y^n_{\delta+\zeta}$

Let $E_{n+1}$ be a club of $\theta$ included in $E_n$ and in each
$E^n_\zeta$ and choose $\overline y^{n+1}$ such that: its range include the
range of $\overline y^n$ and for $\delta_1 <\delta_2 \in E_{n+1}$ we have
$\{ y^{n+1}_{\delta+\xi } \} : \xi <\sigma\}$ include
$\{y^n_{\delta+\xi}:\xi<\sigma\} \cup Z^n_\delta$ and $\{
y^{n+1}_i:i<\delta\}$ include each $z^n_{\zeta, \epsilon}$ for
$\zeta<\sigma$ and $\epsilon<\delta$ .  The rest as as in the proof of
2.3(2).
\newpage

REFERENCES
\bigskip

\bibliographystyle{lit-plain}
\bibliography{lista,listb,listx}
\shlhetal
      
\enddocument
\bye